


\overfullrule=0pt

%
%
%
%
%
%
%
%

%
%
 
 
\baselineskip=14pt     
\parskip=6pt           
\parindent=20pt        


 
\def\today{\ifcase\month\or January\or February\or March\or April\or May\or 
June\or July\or August\or September\or October\or November\or December\fi 
\space\number\day, \number\year}


\font\medbold=cmbx10 scaled 1300

\def\cheading#1{  \centerline{{\bf #1}}   } 
 
\def\subheading#1{\goodbreak\medskip\noindent{\bf #1}\nobreak}
  
\def\section#1 
   {\goodbreak\medskip 
   {\noindent\hglue-1.5pt\newsection\ \bf#1} 
   \nobreak 
} 
 
 
\newcount\aREGISTER  
 
\def\newsection{ 
       \global\advance\aREGISTER by1        
       \global\bREGISTER=0                   
       \global\equationREGISTER=0            
       {${{\bf\the\aREGISTER}}.\ $}          
       } 
                     
\def\sectionno{\the\aREGISTER}  
 
 
\def\list{ 
      \vglue8pt 
      \begingroup 
      \parindent=1cm 
      \parskip=3pt 
      \baselineskip=10pt   
      \resetitemregister 
      \vglue0pt plus5pt 
      \interlinepenalty=20 
      \def\i{\item{$\newitem $.\ }}   
      \def\j{\item{$\newitem $)\ }}}  
    
\def\endlist{ 
      \par\endgroup\goodbreak\vglue0pt plus5pt minus1pt} 
       
\newcount\itemregister  
 
\def\resetitemregister{\itemregister=0} 
 
\def\newitem{ 
      {\global\advance\itemregister by1} 
      \the\itemregister}                 

\def\itemno{\the\itemregister}
 
 
\def\proclaim#1 
  {\medbreak 
   \smallskip 
   \noindent 
   {\bf \newproc.}\ {\bf#1}\   
   \begingroup\sl\penalty80} 
 
\def\endproclaim{\endgroup\smallskip\goodbreak} 
 
\def\pproclaim#1   
  {\medbreak 
   \smallskip 
   \noindent 
   {\bf#1}\   
   \begingroup\sl\penalty80}

%
 
%

\newcount\bREGISTER 
 
\def\proclabel{{\the\aREGISTER}.{\the\bREGISTER}} 
 
\def\newproc{\global\advance\bREGISTER by 1 
     \proclabel 
                     }

 
\def\eqn{\eqno {(\newequation)} } 
 
\newcount\equationREGISTER 
 
\def\eqlabel{{\the\aREGISTER}.{\the\equationREGISTER}}
 
\def\newequation{ 
    \global\advance\equationREGISTER by 1 
    \eqlabel
}


\def\nfootnote#1{ 
    \begingroup 
    \parindent=10pt 
    \baselineskip=9pt 
    \global\advance\footnoteregister by 1 
    \footnote{~$^{\the\footnoteregister}$}{#1}%
\endgroup}%
     
\newcount\footnoteregister 
 

 
\def\comment#1{}

%
%

\def\ifundefined#1{\expandafter\ifx\csname #1\endcsname\relax}

\newwrite\LOWTEXauxfile
\newwrite\LOWTEXerrormsg

\def\begindoc{%
   \openin\LOWTEXauxfile \jobname-aux
   \ifeof\LOWTEXauxfile
      \immediate\write\LOWTEXerrormsg{No auxiliary file yet.}%
   \else
      \closein\LOWTEXauxfile\relax
      \input \jobname-aux
   \fi
   \immediate\openout\LOWTEXauxfile=\jobname-aux
   \immediate\write\LOWTEXauxfile{\relax}%
}

\def\enddoc{%
   \closeout\LOWTEXauxfile
}

\def\procref#1{%
   \xdef\LOWTEXwriteref{%
      \write\LOWTEXauxfile{%
         \def \string\LOWTEXREF #1 {\proclabel}%
      }%
   }%
   \LOWTEXwriteref
}

\def\eqref#1{%
   \xdef\LOWTEXwriteref{%
      \write\LOWTEXauxfile{%
         \def \string\LOWTEXREF #1 {\eqlabel}%
      }%
   }%
   \LOWTEXwriteref
}

\def\xref#1{%
   \ifundefined{LOWTEXREF#1}
      \immediate\write\LOWTEXerrormsg{Undefined reference #1.}%
   \else
      \csname LOWTEXREF#1\endcsname
   \fi
}


\def\bibliography{
  \vfill\eject
  \centerline{{\medbold Bibliography }}
  \medskip
  \begingroup
  \parindent0.5cm
  }

\def\bibref#1{%
\xdef\LOWTEXwriteref{%
\write\LOWTEXauxfile{%
\def \string\LOWTEXREF #1 {\bibno}%
}%
}%
\LOWTEXwriteref
}

\newcount\bibREG  

\def\newbibno{%
\global\advance\bibREG by1
\the\bibREG%
}

\def\bibno{\the\bibREG}

\def\bi#1{%
\item{[\newbibno]}\kern-3pt%
\bibref{#1}%
}

\def\endbibliography{
  \endgroup
  }
  
\def\cite#1{[\xref{#1}]}

%
%

\def\picture#1by#2(#3){ 
\vbox to #2 { 
  \hrule width #1 height 0pt depth 0pt \vfill \special{picture #3}} 
}

\def\scaledpicture#1by#2(#3scaled#4){{ 
\dimen0=#1  \dimen1=#2 
\divide\dimen0 by 1000 \multiply\dimen0 by #4 
\divide\dimen1 by 1000 \multiply\dimen1 by #4 
\picture \dimen0 by \dimen1 (#3 scaled #4)}} 
 
\def\cfigure#1by#2(#3scaled#4offset#5:#6) 
  {\medskip 
   \vglue 2mm minus 2mm 
   $$ 
     \hbox{ 
       \hglue#5 
       {\scaledpicture #1 by #2 (#3 scaled #4)} 
     } 
   $$ 
   \par\nobreak 
   \vglue 4.5mm\nobreak 
   \centerline {\figurelabel \quad #6}\par\goodbreak 
   \vglue 2mm minus 2mm 
   \medskip} 
 
\newcount\figureregister 
\def\figurenumber 
  {\global\advance\figureregister by 1 
    \the\figureregister} 
 
\def\figurelabel 
  {{\bf Figure \figurenumber .}} 

%
%
%
%



\def\lset{\{\ }  
\def\st{\ |\ }   
\def\rset{\ \}}  
 
\def\set#1{\lset #1 \rset} 
\def\sett#1#2{\lset #1 \st #2 \rset} 
 
\def\map{\longrightarrow}

\def\mapright#1{\ \smash{ 
   \mathop{\longrightarrow}\limits^{#1}}\ }


\def\Hom{\mathop{\rm Hom}\nolimits}

 
\def\dim{\mathop{\rm dim\, }\nolimits}

 
\def\tr{{}^t\kern-0.9pt} 
 
\def\Tr{{}^T\kern-0.9pt} 
 
\def\wedges{\wedge \dots \wedge}


\def\approxle{ 
       \mathop{  
       \hbox{$ 
       {  \lower0.5ex\hbox{$<$}  
              \atop  
           \raise0.25ex\hbox{$\sim$}   
       } 
             $} 
               } 
              } 
 
\def\del{\partial} 
 
\def\exp{\mathop{\rm exp}\nolimits}


 
\def\commadots{\ ,\ldots ,\ }

 

 
\def\rank{\mathop{{\rm rank}}\nolimits} 
 
\def\perp{^\bot}

\def\numeq{
\vcenter{
\hbox{
\vbox{
  \hbox{$\approx$}%
  \kern-8.6pt%
  \hbox{$\sim$}%
     }
     }
       }
      \,\,%
} 



\def\proof{\noindent{\bf Proof:}}

%
%
%
%
%
%
%
 
\newfam\meuffam 
\font\tenmeuf=eufm10 
\font\sevenmeuf=eufm7 
\font\fivemeuf=eufm5 
 
\textfont\meuffam=\tenmeuf 
\scriptfont\meuffam=\sevenmeuf 
\scriptscriptfont\meuffam=\fivemeuf 
 
\def\germ{\fam\meuffam\tenmeuf}

 

\newfam\msyfam
\font\tenmsy=msbm10
\font\sevenmsy=msbm7
\font\fivemsy=msbm5

\textfont\msyfam=\tenmsy
\scriptfont\msyfam=\sevenmsy
\scriptscriptfont\msyfam=\fivemsy

\def\Bbb{\fam\msyfam\tenmsy}


\def\g{{\germ g}} 
 
\def\i{{\germ i}} 
\def\j{{\germ j}}


\def\C{{\Bbb C}}

\def\H{{\Bbb H}}

\def\P{{\Bbb P}} 
\def\Q{{\Bbb Q}} 
\def\R{{\Bbb R}}

\def\Z{{\Bbb Z}} 
 

\def\LL{{\cal L}}

\def\VV{{\cal V}}



\def\del{\partial}
\def\proof{\noindent{\bf Proof:\ }}
\def\res{\hbox{res}\,}
\def\hide#1{}

\def\Aut{\hbox{Aut}}

\def\bX{{\bf X}}
\def\bY{{\bf Y}}

\def\mod{\hbox{ mod }}

\def\eqrefer#1{$(\xref{#1})$}

\baselineskip=12pt

\font\bigfont=cmbx10 scaled\magstep1

\def\comp{\circ}

\def\bibliography{
  \medskip
  \subheading{{\bf References }}
  \smallskip
  \begingroup
  \parindent1cm
  }

\def\bi#1{%
\item{[\newbibno]}\kern-3pt%
\bibref{#1}%
}

\def\secref#1{%
   \xdef\LOWTEXwriteref{%
      \write\LOWTEXauxfile{%
         \def \string\LOWTEXREF #1 {\sectionno}%
      }%
   }%
   \LOWTEXwriteref
}


\def\cokernel{\mathop{cokernel}}
\def\kernel{\mathop{kernel}}

\parskip=0pt


\begindoc

\cheading{{ \bigfont Discriminant Complements and }}
\cheading{{ \bigfont Kernels of Monodromy Representations }}

\medskip
\cheading{James A. Carlson and Domingo Toledo}
\footnote{}{Authors
partially supported by National Science Foundation Grant DMS 9625463}


\medskip
\begingroup
\narrower\narrower

\cheading{Abstract}

Let $\Phi_{d,n}$ be the fundamental group of the space of smooth
projective hypersurfaces of degree $d$ and dimension $n$ and let
$\rho$ be its natural monodromy representation.  Then the kernel
of $\rho$ is {\sl large} for $d \ge 3$ with the exception of 
the cases $(d,n) = (3,0),\ (3,1)$.  For these and for  $d < 3$
the kernel is finite.  A large group is one that admits a
homomorphism to a semisimple Lie group of noncompact type with Zariski-dense
image.  By the Tits alternative a large group contains a free subgroup of rank
two.  

\endgroup

\section{Introduction}
\secref{introsection}

A hypersurface of degree $d$ in a complex projective space
$\P^{n+1}$ is defined by an equation of the form
$$
    F(x) = \sum a_L x^L = 0,
    \eqn
    \eqref{universalhypersurface}
$$
where $x^L = x_0^{L_0} \cdots x_{n+1}^{L_{n+1}}$ is a monomial of degree
$d$ and where the $a_L$ are arbitrary complex numbers, not all zero. Viewed as
an equation in both the $a$'s and the $x$'s, \eqrefer{universalhypersurface} 
defines a hypersurface $\bX$ in $\P^N\times\P^{n+1}$, where $N+1$ is the
dimension of the space of homogeneous polynomials of degree $d$ in $n+2$
variables, and where the projection
$p$ onto the first factor makes $\bX$ into a family with fibers
$X_a = p^{-1}(a)$. This is the universal family of hypersurfaces of degree $d$
and dimension
$n$.  Let $\Delta$ be the set of points $a$ in $\P^N$ such that the
corresponding fiber is singular.  This is the {\sl discriminant locus}; it
is well-known to be irreducible and of codimension one. Our aim is to study
the fundamental group of its complement, which we write as
$$
  \Phi = \pi_1(\P^N - \Delta).
$$
When we need to make precise statements we will sometimes write
$
  \Phi_{d,n} = \pi_1(U_{d,n}, o),
$
where $d$ and $n$ are as above, $U_{d,n} = \P^N - \Delta$, and $o$ is a base
point.

The groups $\Phi$ are almost always  nontrivial and in fact are almost
always {\sl large}.  By this we mean that there is a homomorphism of $\Phi$
to  a non-compact semi-simple real algebraic group which has Zariski-dense
image.  Large groups are infinite, and, moreover, always contain a free group of
rank two.  This follows from the Tits alternative
\cite{Tits}, which states that in characteristic zero a linear group  either
has a solvable subgroup of finite index or contains a free group of rank
two.

To show that $\Phi = \Phi_{d,n}$ is large we consider the image 
$\Gamma = \Gamma_{d,n}$ of the
monodromy representation
$$
  \rho: \Phi \map G.
\eqn
\eqref{monodrep}
$$
Here and throughout this paper  $G = G_{d,n}$ denotes the group of
automorphisms of the primitive cohomology $H^n(X_o,\R)_o$ which preserve the
cup product. When $n$ is odd the primitive cohomology is the
same as the cohomology, and when $n$ is even it is the orthogonal complement
of $h^{n/2}$,  where $h$  is the hyperplane class.  Thus $G$ is either a
symplectic or an orthogonal group, depending on the parity of $n$, and is an 
almost simple  real  algebraic group. 

About the image of the monodromy representation, much is known.  Using
results of Ebeling \cite{Ebeling} and Janssen \cite{Janssen}, Beauville in \cite{BeauvilleLattice}
established the following:

\proclaim{Theorem}.  Let $G_\Z$ be the subgroup of $G$ which preserves the
integral cohomology.  Then the monodromy group $\Gamma_{d,n}$ is of finite index in
$G_\Z$.  Thus it is an arithmetic subgroup.
\endproclaim
\procref{Beauvillethm}

\noindent 
The result in \cite{BeauvilleLattice} is much more
precise: it identifies $\Gamma$ as a specific subgroup of finite
(and small) index in $G_\Z$. Now suppose that $d > 2$ and that $(d,n) \ne (3,2)$. 
Then $G$ is noncompact, and the results of Borel  \cite{BorelDT} and  
Borel-Harish-Chandra \cite{BHC}
apply to show that
$\Gamma$ is (a) Zariski-dense and (b) a lattice.  Thus (a) the smallest algebraic
subgroup of $G$ which contains $\Gamma$ is $G$ itself and (b) $G/\Gamma$ has
finite volume.

Consider now the kernel of the monodromy representation, which we denote by
$K$ and which fits in the exact sequence
$$
   1 \map K \map \Phi \mapright{\rho} \Gamma \map 1.
   \eqn
   \eqref{KPhiGamma}
$$
The purpose of this paper is to show that in almost
all cases it is also large:

\proclaim{Theorem.} The kernel of the monodromy
representation \eqrefer{monodrep} is large if $d > 2$ and $(d,n) \ne (3,1),\
(3,0)$.
\endproclaim
\procref{maintheorem}

The theorem is sharp in the sense that the remaining groups are finite.  
When $d = 2$, the case of  quadrics, $\Phi$ is finite cyclic.  
When $(d,n) = (3,0)$,  the configuration space $U$ parametrizes
unordered sets of three distinct points in the projective line and so $\Phi$
is the braid group for three strands in the sphere.  It has order
12 and can be faithfully represented by
symmetries of  a regular hexagon. 

When $(d,n) = (3,1)$ the configuration space $U$ parametrizes smooth cubic
plane curves and the above sequence can be written as
$$
   1 \map K \map \Phi_{3,1} \mapright{\rho} SL(2,\Z) \map 1,
$$
where $K$ is the three-dimensional Heisenberg group over the field $\Z/3$, a
finite group of order 27. Moreover, $\Phi_{3,1}$ is a semi-direct product,
where $SL(2,\Z )$ acts on $K$ in the natural way. This result, due
to Dolgachev and Libgober \cite{DolgLib}, is to our knowledge the only one 
which determines
the exact sequence \eqrefer{KPhiGamma} for hypersurfaces of positive
dimension and degree larger than two.   
Note that in this case $\Phi$ is large but $K$ is finite.

Note also that there are two kinds of groups for which the natural monodromy
representation has finite image but large kernel.  These are the braid
groups $\Phi_{d,0}$ for $d > 3$ and the group $\Phi_{3,2}$ for the space of
cubic surfaces. Thus all of them are large.  For the braid groups
this result is classical, but for $\Phi_{3,2}$ it is new.  Since $\Phi_{3,2}$
is large  it is infinite, a fact which answers a question left open
by Libgober in \cite{Lib}.  

Concerning the proof of Theorem \xref{maintheorem}, we would like to say first of
all that it depends, like anything else in this subject, on the Picard-Lefschetz
formulas.  We illustrate their importance by sketching how they imply
the non-triviality of the monodromy representation
\eqrefer{monodrep}.  Consider a smooth point
$c$ of the discriminant locus.  For these $X_c$ has a
exactly one node: an isolated singularity defined in suitable local coordinates
by a nondegenerate sum of squares.  Consider also a loop $\gamma = \gamma_c$
defined by following a path $\alpha$ from the base point to the edge of a
complex disk normal to $\Delta$ and centered at $c$, traveling once around the
circle bounding this disk, and then returning to the base point along
$\alpha$ reversed.  By analogy with the case of knots, we
call these loops (and also their homotopy classes) the {\sl meridians} of
$\Delta$.  Then $T = \rho(\gamma)$ is a {\sl Picard-Lefschetz} transformation,
given by the formula
$$
   T(x) = x \pm (x,\delta) \delta .
   \eqn
   \eqref{plformula}
$$
Here $(x,y)$ is the cup product and $\delta$ is the {\sl vanishing cycle}
associated to $\gamma$.  When $n$ is odd, $(\delta,\delta) = 0$ and the sign
in \eqrefer{plformula} is $-$.  When $n$ is even and $(\delta,\delta) = \pm
2$, the sign in \eqrefer{plformula} is $\mp$ (see \cite{DeWeOne}, paragraph
4.1). Thus when $n$ is even $\delta$ is automatically nonhomologous to zero, and so
$T$ must be nontrivial.  Since vanishing cycles exist  whenever the hypersurface
$X_o$ can degenerate to a variety with a node, we conclude that $\rho$ is
nontrivial for $n$ even and $d > 1$.  Slightly less elementary arguments show that
the homology class of the vanishing cycle, and hence the monodromy representation,
is nontrivial for all $d > 1$ except for the case $(d,n) = (2,1)$.

The proofs of theorem \xref{Beauvillethm}, an earlier result of
Deligne asserting the Zariski density of $\Gamma_{d,n}$, and the main result of
this paper are based on the Picard-Lefschetz formulas
\eqrefer{plformula}.  Our proof begins with the construction of a universal family
of cyclic  covers of $\P^{n+1}$ branched along the hypersurfaces
$X$.  From it we define a second monodromy representation $\bar \rho'$ of
$\Phi$.  Suitable versions of the Picard-Lefschetz formulas and
Deligne's theorem apply to show that $\bar \rho'$ has Zariski-dense
image.  Finally, we apply Margulis' super-rigidity theorem to show that
$\bar\rho'(K)$, where $K$ is the kernel of the natural monodromy representation, is
Zariski-dense.  Thus $K$ is large.

We mention the paper 
\cite{Mag} as an example of the use of an associated family of 
cyclic covers  to
construct  representations (in this case for the braid
groups of the sphere).  We also  note the related results of the
article \cite{DOZ} which we learned of while preparing the final version of this
manuscript.  The main theorem is that the complement of the dual $\widehat
C$ of an immersed curve $C$ of genus at least one, or of an immersed rational
curve of degee at least four, is {\sl big} in the sense that it contains a free
group of rank two.  When $C$ is smooth, imbedded, and of even degree at least
four this follows from a construction of Griffiths \cite{GriffHyperbolic}:
consider the family of hyperelliptic curves obtained as double covers of a line
$L$ not tangent to $C$ which is branched at the points $L\cap C$.  It defines
a monodromy representation  of $\Phi = \pi_1(\widehat{\P}^2 - \widehat{C})$ with
Zariski-dense image.  Consequently $\Phi$ is large, and, {\sl a fortiori},
big. Such constructions have inspired the present paper. By using cyclic 
covers of higher degree one can treat the case of odd degree greater than
four in the same way.

The authors would like to thank Herb Clemens and Carlos Simpson for
very helpful discussions.

\section{Outline of the proof}
\secref{outlinesection}

As noted above, the proof of the main theorem is based on the construction of an
auxiliary representation $\rho'$ defined via a family of cyclic covers $Y$ of
$\P^{n+1}$ branched along the hypersurfaces $X$.  To describe it, let
$k$ be a divisor of $d$ and consider the equation
$$
    F(a,x) = y^k + \sum a_L x^L = 0,
    \eqn
    \eqref{universalcyclic}
$$
which for the moment we view as defining
a set $\widehat \bY$ in $(\C^{N+1} - \set{0})\times \C^{n+3}$ with coordinates
$a_L$ for $\C^{N+1}$ and coordinates $x_0,\cdots ,x_{n+2}$ and $y$ for
$\C^{n+3}$. Construct an action of
$\C^*$ on it by  multiplying the coordinates $x_i$ by $t$ and by multiplying
$y$ by $t^{d/k}$.  View the quotient $\bY$ in
$(\C^{N+1} - \set{0})\times \P^{n+2}$, where we use $\P^{n+2}$ to denote
the weighted projective space for which the $x_i$ have weight one and for which $y$
has weight $d/k$.  

The resulting universal family of cyclic covers $\bY$ is defined on
$\C^{N+1} - \set{0}$ and has smooth fibers over $\widetilde U =  \C^{N+1} -
\widetilde
\Delta$, where $\widetilde \Delta$ is the pre-image of $\Delta$. Since $\C^{N+1} -
\set{0}$ is a principal
$\C^*$ bundle over $\P^N$, the same holds over  $\widetilde U$ and $\widetilde
\Delta$. It follows that one has a central extension
$$
  0 \map \Z \map \widetilde \Phi \map \Phi \map 1,
$$  
where $\widetilde \Phi = \pi_1(\widetilde U)$.  We introduce
$\widetilde U$ and $\widetilde\Phi$ purely for the technical reason that the
universal family  of cyclic branched covers need not be defined over $U$ itself.

The family $\bY |\widetilde U$ has a monodromy representation which
we denote by
$\tilde\rho$ and which takes values in a real algebraic group $\widetilde G$ of
automorphisms of $H^{n+1}(Y_{\tilde o},\C )$ which commute with the cyclic
group of covering transformations (and which preserve the hyperplane class and
the cup product).  Here $\tilde o$ is a base point in $\widetilde U$ which lies
above the previously chosen base point $o$ of $U$, and $Y_{\tilde o}$ denotes the
$k$-fold cyclic cover of $\P^{n+1}$ branched over $X_o$.  

The group $\widetilde G$ is semisimple but in general has more than one
simple factor.  Let $G'$ be one of these and let
$$
\rho ' :\widetilde\Phi\map G',
$$
denote the  composition of $\tilde\rho$ with the projection to $G'$.  Then
we must establish the following:

\proclaim{Technical point}. The factor $G'$ can be chosen to be a non-compact
almost simple real algebraic group.  The image of
$\rho '$ is Zariski-dense in
$G'$.
\endproclaim
\procref{technicalpoint}

 Suppose that this is true.  Then we can argue as follows.  First, the group
of matrices which commute with
$\rho'(\Z)$ contains a Zariski-dense group.  Consequently $\rho'(\Z)$ lies in the
center of
$G'$.  Therefore there is a quotient representation
$$
  \bar \rho': \Phi\map \bar G',
$$
where $\bar G' $ is the adjoint group of $G'$ (that is, $G'$ modulo its center).  
Moreover, the  representation
$\bar\rho'$ also has Zariski-dense image.  

Now consider our original representation \eqrefer{monodrep}.    Replacing
$\Phi$ by a normal subgroup of finite index we may assume that the image of
$\rho$ lies in the identity component of $G$ in the analytic topology and
that the image of $\bar\rho '$ lies in the identity component of $\bar G'$
in the Zariski topology.  Let $\bar G$ denote the identity component (in the
analytic topology) of $G$ modulo its center, and   let
 $\bar\rho :\Phi\map \bar
G$ denote the resulting representation.   We still have that
$\bar\rho (\Phi)$ is a lattice in
$\bar G$ and that $\bar \rho '(\Phi)$ is Zariski-dense in $\bar G'$.

Now let $\bar K$ be the kernel of $\bar\rho$, and let $L$ be the
Zariski-closure of $\bar\rho'(\bar K)$.  Since $\bar K$ is normal in $\Phi$
and $\bar\rho'(\Phi)$ is Zariski-dense in $\bar G'$, $L$ is normal in $\bar
G'$.  Since $\bar G'$ is  a {\sl simple} algebraic group, either $L = \bar G'$
or $L = \set{1}$.  If the first of the two alternatives holds, then 
$\bar\rho'(K)$ is Zariski dense, and so $K$ is large. This is because $K$ has
finite index in $\bar K$ and so $\bar\rho'(K)$ and $\bar\rho'(\bar K)$ have
the same Zariski closure.  

We now show that the second alternative leads to a contradiction, from which it
follows that $K$ must be large.  Indeed, if $\bar \rho'(\bar K) =
\set{1}$, then the expression $\bar \rho'\comp\bar \rho^{-1}$
defines a homomorphism from the lattice $\bar \rho(\Phi)$ in $\bar G$ to the
Zariski-dense subgroup $\bar \rho'(\Phi)$ in
$\bar G'$.   If the real rank of $\bar G$ is at least two, the Margulis
rigidity theorem
\cite{MargulisRigidity}, \cite{Zimmer} Theorem 5.1.2, applies to give an
extension of $\bar \rho'$ to a homomorphism of $\bar G$ to
$\bar G'$.  Since $\bar\rho'(\Phi)$ is Zariski-dense, the extension is surjective. 
Since $\bar G$ is simple, it is an isomorphism. Thus the complexified lie algebras
$\g_\C, \g'_\C$ must be isomorphic.   However, one easily shows that $\g_\C
\not\cong \g'_\C$, and this contradiction completes the proof.

We carry out the details separately in two cases.  First, for the simpler case
where $d$ is even and $(d,n)\ne (4,1)$, we use double covers ($k = 2$).  Then
$G'$ is the full group of automorphisms of the
primitive (or anti-invariant) part of  $H^{n+1}(Y,\R)$ and so is again an
orthogonal or symplectic group.  The technical point \xref{technicalpoint}
follows from a density result of Deligne that we recall in section
\xref{zardensitysection}.  Deligne's result gives an alternative between
Zariski 
density and finite image, and the possibility of finite image is excluded in
section
\xref{ratdiffsection}.  Finally the  Lie algebras $\g_\C$ and
$\g'_\C$ are not isomorphic, since when one of them is symplectic (type
$C_\ell$), the other is orthogonal (type $B_\ell$ or $D_\ell$).  By lemma
\xref{rankboundslemma} the rank $\ell$ is at least three, so there are no
accidental isomorphisms, e.g., $B_2 \cong C_2$.

For the remaining cases, namely $d$ odd or $(d,n) = (4,1)$ we use $d$-fold
covers, i.e., $k = d$.    For these we must identify the group $\widetilde G$
of automorphisms of $H^{n+1} (Y,\R)_0$ which preserve the cup product and
which commute with the cyclic automorphism $\sigma$.  This is the natural 
group in which the monodromy representation $\tilde \rho$ takes its values.
Now a linear map commutes with $\sigma$ if and only if it preserves the
eigenspace decomposition of $\sigma$, which we write as
$$
 H^{n+1} (Y,\C)_0 = \bigoplus_{\mu \ne 1} H(\mu).
$$
As noted in \eqrefer{eigenspacedimform}, the dimension of $H(\mu)$
is independent of $\mu$.  Now let $\widetilde G(\mu)$ be subgroup of $\widetilde G$ which acts by
the identity on $H(\lambda)$ for $\lambda \ne \mu,\; \bar\mu$.  It can be viewed as a
group of transformations of $H(\mu) + H(\bar\mu)$.  Thus there is a
decomposition
$$
   \widetilde G 
    = 
     \prod_{ \mu \in S} \widetilde G(\mu),
   \eqn
   \eqref{tildeGdecomp}
$$
where
$$
   S = \sett{ \mu }{ \mu^k =1,\ \mu \ne 1,\ \Im \mu \ge 0 }.
$$
When $\mu$ is non-real, $\widetilde G(\mu)$ can be identified 
via the projection $H(\mu) \oplus H(\bar\mu) \map H(\mu)$
with the group of transformations of $H(\mu)$ which are unitary with respect
to the hermitian form $h(x,y) = i^{n+1}(x,\bar y)$, where $(x,y)$ is the cup
product.  This form may be (and usually is) indefinite.  When
$\mu = -1$, $\widetilde G(\mu)$ is the group of transformations of $H(-1)$
which preserve the cup product.  It is therefore an orthogonal or symplectic group.

We will show that at least one of the components $\widetilde\rho_\mu(\Phi)\subset
\widetilde G(\mu)$ is Zariski-dense, and we will take $G' = \widetilde G(\mu)$.  
The necessary Zariski density result, which is a straightforward adaptation of
Deligne's, is proved in section \xref{unitarydensitysection} after some preliminary
work on complex reflections in section \xref{complexreflectionsection}.  Again,
the possibility of finite image has to be excluded, and the argument for this is
in section \xref{cycliccoversection}.  Finally, to prove that
$\g_\C$ and $\g'_\C$ are not isomorphic one observes that
$\g_\C$ is of type $B_\ell, C_\ell$ or $D_\ell$ while $\g'_\C$ is of type
$A_\ell$ (since
$G'$ is of type $SU(r,s)$.  One only needs to avoid the isomorphism $D_3\cong
A_3$, which follows from the lower bound of the rank of $\g_\C$ in lemma
\xref{rankboundslemma}.

In order to apply Margulis' theorem we also need to verify that the real
rank of $G$ is at least two.  This is done in section \xref{rankbounds}.   

To summarize, we have established the following general criterion, and our proof of
Theorem \xref{maintheorem} is an application of it.

\proclaim{Criterion.}  The kernel $K$ of a linear representation $\rho: \Phi \map
G$ is large if
\list
\i $\rho(\Phi) \subset G$ is a lattice in a simple Lie group $G$ of
real rank at least two.
\i There exist a non-compact, almost simple real algebraic group $G'$, a central 
extension $\widetilde\Phi$ of $\Phi$ and a linear 
representation
$\rho':
\widetilde
\Phi
\map G'$  with Zariski-dense image.
\i $G$ and $G'$ are not locally isomorphic.
\endlist
\endproclaim
\procref{largekernelcriterion}

\noindent
An immediate consequence is the following:

\proclaim{Corollary.} Let $\Phi$ be a group which admits a representation 
$\rho: \Phi \map G$
to a simple Lie group of real rank greater than 1 with image a lattice.  
Suppose further that there exist an almost simple real algebraic group $G'$, 
 a central extension $\widetilde \Phi$ of $\Phi$, and a representation $\rho':
\widetilde\Phi \map G'$ with  Zariski-dense image.
Suppose in addition that $G$ and $G'$ are not locally isomorphic.  Then $\Phi$ is not
isomorphic to a lattice in any simple Lie group of real rank greater than 1.
\endproclaim

\proof Suppose that $\tau: \Phi \map \Sigma$ is an isomorphism of $\Phi$
with a lattice $\Sigma$ in a Lie group $H$ of real rank greater than one.
If $H$ is not locally isomorphic to $G'$, then apply the criterion with $\tau$ in 
place of $\rho$ to conclude that $\tau$ has large kernel, hence cannot be an
isomorphism.  Suppose next that $H$ is locally isomorphic to $G'$.  Apply 
the criterion with $\tau$ in place of $\rho$ and with $\rho$ in place of $\rho'$
to conclude as before that the kernel of $\tau$ is large.

For most families of hypersurfaces the natural monodromy representation and the
representation for the associated family of cyclic covers satisfy the hypotheses of
the corollary to give the following:

\proclaim{Theorem.} If $d > 2$, $n>0$, and $(d,n) \ne (3,1),\; (3,2)$, the group
$\Phi_{d,n}$ is not isomorphic to a lattice in a simple Lie group of real rank
greater than one.
\endproclaim

It seems reasonable that the preceding theorem holds with ``semisimple''
in place of ``simple.''  However, we are unable show that this is the case.
Indeed, our results so far are compatible with an isomorphism $\Phi \cong
\Gamma\times\Gamma'$.  We can exclude this in certain cases (see section
\xref{remarkssection}), but not for an arbitrary subgroup of finite index, which
is what one expects.

\proclaim{Remarks.}
\procref{remarkMeridiansInfiniteOrder}
\rm

(a) Suppose that $d \ge 3$ and let $\gamma$ be a meridian of $\Phi_{d,n}$. 
When $n$ is odd, $\rho(\gamma)$  is a nontrivial symplectic transvection. Since
it is of infinite order, so is the  meridian $\gamma$.  When $n$ is even,
$\rho(\gamma)$ is a reflection, hence of order two.  Now suppose that $d$ 
is even and consider the monodromy representation of the central extension $\widetilde
\Phi$ constructed from double covers. Let $\tilde\gamma$ be a lift of 
$\gamma$ to an element of $\widetilde\Phi$. Then
$\rho'(\tilde\gamma)$ is a nontrivial symplectic transvection, no power of which is central. Thus
$\bar\rho'(\gamma)$ is of infinite order, and, once again, we conclude that
$\gamma$ is of infinite order.

(b) M. Kontsevich informs us that he can prove that for any $d>2$ (and at
least for $n=2$) the local monodromy corresponding to a meridian is  of
infinite order in the group of connected components of the symplectomorphism
group of $X_o$.  This implies that the meridians are of infinite order for all
$d>2$, not necessarily even as above. The symplectic nature of the monodromy
for a meridian (for $n=2$) is studied in great detail by P. Seidel in his
thesis \cite{Seidel}.

(c) For the case of double covers the image $\Gamma'$ of the fundamental group
under the second monodromy representation $\rho'(\widetilde\Phi)$ is a
lattice. This follows from the argument given by Beauville to prove theorem
\xref{Beauvillethm}. It is enough to be able to degenerate the branch locus
$X$ to a variety which has an isolated singularity of the form $x^3 + y^3 + z^4
+ \hbox{a sum of squares} = 0$. Then the roles of the kernels $K$ and $K'$ are
symmetric and  one concludes that $K'$ is also large.

\endproclaim

\section{Zariski Density}
\secref{zardensitysection}

The question of Zariski-density for monodromy groups of Lefschetz
pencils was settled by Deligne in \cite{DeWeOne} and \cite{DeWeTwo}.
We review these results here in a form convenient for the proof of the main
theorem in the case of even degree and also for the proof of a
density theorem for unitary groups (section
\xref{unitarydensitysection}).  To begin, we have the following
purely group-theoretic fact:  \cite{DeWeTwo}(4.4):

\proclaim{Theorem. (Deligne)}  Let $V$ be a vector space (over $\C$)
with a non-degenerate bilinear form $(\ ,\ )$ which is either symmetric or
skew-symmetric.  Let
$\Gamma$ be a group of  linear transformations of $V$ which preserves the
bilinear form.   Assume the existence of a subset $E\subset V$ such that
$\Gamma$ is generated by the Picard-Lefschetz transformations \eqrefer{plformula}
with $\delta\in E$.  
  Suppose that $E$ consists of a single $\Gamma$-orbit
and that it spans $V$.  Then $\Gamma$ is either finite or Zariski-dense.
\endproclaim
\procref{delignedensity}

To apply this theorem in a geometric setting, consider a family of
$n$-dimensional varieties $p: \bX \map S$  with discriminant locus $\Delta$ and 
monodromy representation
$\rho:\pi_1(S - \Delta) \map \Aut(H^n(X_o))$.  Assume that $S$ is either
$\C^{N+1} - \{ 0\} , N\ge 1$ or $\P^N$, so that $S$ is simply connected and 
hence that
$\pi_1 (S-\Delta)$ is generated by {\sl meridians} (cf. \S
\xref{introsection} for the definition).  Assume also that for
each meridian there is a class $\delta\in H^n (X_o )$ such that the
corresponding monodromy transformation is given  by the Picard-Lefschetz
formula 
\eqrefer{plformula}.  Let $E$ denote the set of these classes (called the 
{\sl vanishing cycles}). Let $V^n(X_o)\subset H^n (X_o)$ be the span of 
$E$, called the {\sl vanishing cohomology}.

A cycle orthogonal to $V = V^n(X_o)$ is invariant under all Picard-Lefschetz
transformations, hence is invariant under the action of monodromy.
Consequently its orthogonal complement $V\perp$ is the space of invariant 
cycles.  The image of $H^n(\bX)$ in $H^n(X_o)$ also consists
of invariant cycles.  By theorem 4.1.1 (or corollary (4.1.2)) of
\cite{DeHodgeTwo}, this inclusion is an equality.  One concludes that $V\perp$
is the same as the image of $H^n(\bX)$, which is a sub-Hodge structure, and so
the bilinear form restricted to it is nondegenerate.  Therefore the
bilinear form restricted to $V = V^n(X_o)$ is also nondegenerate.  Consequently 
$V^n(X_o)$ is an orthogonal or symplectic space, and the monodromy group acts
on $V^n(X_o)$ by orthogonal or symplectic transformations.

When the discriminant locus is irreducible the argument of Zariski
\cite{Zar} or \cite{DeWeOne}, paragraph preceding Corollary 5.5, shows that the
meridians of $\pi_1(S - \Delta)$ are mutually conjugate. Writing down a
conjugacy $\gamma' = \kappa^{-1}\gamma\kappa$ and applying it to
\eqrefer{plformula}, one concludes that $\delta' =
\rho(\kappa^{-1})(\delta)$.  Thus the vanishing cycles constitute a single
orbit.  To summarize, we have the following, (c.f.
\cite{DeWeOne}, Proposition 5.3, Theorem 5.4, and \cite{DeWeTwo}, Lemma 4.4.2):

\proclaim{Theorem.} Let $\bX \map S$, with $S = \C^{N+1} - \{ 0\}$ or $\P^N$ 
and $N \ge 1$, be a
family with irreducible discriminant locus and such that the monodromy
transformations of meridians are Picard-Lefschetz transformations.  Then
the monodromy group is  either finite or is a Zariski-dense subgroup of the 
(orthogonal or symplectic) group of automorphisms of the vanishing cohomology.
\endproclaim
 
To decide which of the two alternatives holds, consider the period mapping
$$
   f : U \map D/\Gamma,
$$
where $D$ is the space \cite{GriffPerDom} which classifies the Hodge
structures $V^n(X_a)$ and where $\Gamma$ is the monodromy group. Then one
has the following well-known principle:

\proclaim{Lemma.} If the monodromy group is finite, then the period map is
constant.
\endproclaim
\procref{finiteimagelemma}

\proof Let $f$ be the period map and suppose that the monodromy representation
is finite.  Then there is an unramified cover $\widetilde S$ of the domain of
$f$  for which the monodromy representation is trivial.  Consequently there is
lift
$\tilde f$ to $\widetilde S$ which takes values in the period domain $D$.  Let
$\bar S$ be a smooth compactification of
$\widetilde S$.  Since $D$ acts like a bounded domain for horizontal holomorphic
maps,
$\tilde f$ extends to a holomorphic map of $\bar S$ to $D$.  
Any such map with compact domain is constant \cite{GS}.

As a consequence of the previous lemma and theorem, we have a practical
density criterion:

\proclaim{Theorem.} Let $\bX$ be a family of varieties over $\C^{N+1} - \{ 0\}$
or
$\P^N$, $N\ge 1$, whose monodromy group is generated by Picard-Lefschetz
transformations
 \eqrefer{plformula},
which has irreducible discriminant locus, and whose period map has nonzero
derivative at one point.  Then the monodromy group is Zariski-dense in the
(orthogonal or symplectic) automorphism group of the vanishing cohomology.
\endproclaim
\procref{practicaldensitycriterion}

Irreducibility of the discriminant locus for hypersurfaces is well known, and 
can be proved as follows.  
Consider the Veronese imbedding
$v$ of $\P^{n+1}$ in $\P^N$. This is the map which sends the homogeneous
coordinate vector $[ x_0 \commadots x_{n+1} ]$ to $[ x^{M_0} \commadots x^{M_N}
]$ where the $x^{M_i}$ are an ordered basis for the monomials of degree $d$ in
the $x_i$.  If $H$ is a hyperplane in $\P^N$, then $v^{-1}(H)$ is a hypersurface
of degree $d$ in $\P^{n+1}$.  All hypersurfaces are obtained in this way,
so the dual projective space $\widehat \P^N$ parametrizes the universal
family.  A hypersurface is singular if and only if $H$ is tangent to 
the Veronese manifold $\VV = v(\P^{n+1})$.  Thus the discriminant 
is the variety $\widehat \VV$ dual to $\VV$.  Since the  variety dual to an 
irreducible variety is also irreducible, it follows that the discriminant is 
irreducible.

Finally, we observe that in the situations considered in this paper, vanishing 
cohomology and primitive cohomology coincide.  This can easily be checked by 
computing the invariant cohomology using a suitable compactification and appealing 
to (4.1.1) of \cite{DeHodgeTwo}.  Since this is not essential to our 
arguments we omit further details.

\section{Rational differentials and the Griffiths residue calculus}  
\secref{ratdiffsection}

Griffiths' local Torelli theorem \cite{GriffPerRat} tells us
that the period map for hypersurfaces of degree $d$ and dimension $n$ is
is nontrivial for $d > 2$ and $n > 1$ with the exception of the case $(d,n) = (3,2)$.
In fact, it says more: the kernel of the differential is the tangent space
to the orbit of the natural action of the projective linear group.  The proof 
is based on the residue calculus for rational differential forms and some
simple commutative algebra (Macaulay's theorem).

What we require here is a weak (but sharp) version of Griffiths' result for
the variations of Hodge structures defined by families of cyclic covers of
hypersurfaces.  For double covers this is straightforward, since such covers
can be viewed as hypersurfaces in a weighted projective space \cite{Dolgachev}.
For higher cyclic covers the variations of Hodge structure are complex, and in 
general the symmetry of Hodge numbers, $h^{p,q} = h^{q,p}$ is broken.  Nonetheless,
the residue calculus still gives the needed result.  Since this last part
is nonstandard, we sketch recall the basics of the residue calculus, how
it applies to the case of double covers, and how it extends to the case of
higher cyclic covers.

To begin, consider weighted projective space $\P^{n+1}$ where the weights of
$x_i$ are $w_i$.  Fix a weighted homogeneous polynomial $P(x)$ and let $X$
be the variety which it defines.  We assume that it is smooth.  Now take a
meromorphic differential $\nu$ on $\P^{n+1}$ which has
a pole of order $q+1$ on $X$.
Its residue is the cohomology class on $X$ defined by the formula
$$
   \int_\gamma \res \nu = { 1 \over 2 \pi } \int_{\partial T(\gamma)} \nu,
$$
where $T(\gamma)$ is a tubular neighborhood of an $n$-cycle $\gamma$.  The
integrand can be written as
$$
   \nu(A,P,q) = { A\Omega \over P^{q+1} }. \eqn\eqref{ratdiff}
$$
where
$$
  \Omega = \sum (-1)^i\;w_i x_i \;dx_0 \wedges \widehat{dx_i} \wedges dx_{n+1}.
$$
The ``volume form'' $\Omega$ has weight $w_0 + \cdots + w_{n+1}$ and the degree of $A$,
which we write as $a(q)$, is such that $\nu$ is of weight zero.  The  primitive
cohomology of $X$ is spanned by Poincar\'e residues of rational differentials, 
and the space of residues with a pole of order $q+1$ is precisely
$F^{n-q}H^n_o(X)$, the $(n-q)$-th level of the Hodge filtration on the primitive cohomology.  When
the numerator polynomial is a linear combination of the partial derivatives of $P$, the residue
is cohomologous in $\P^{n+1} - X$ to a differential with a pole of order one
lower.  Let $J = (\del P/\del x_0 \commadots \del P/\del x_{n+1})$
be the Jacobian ideal and let $R = \C[x_0 \commadots x_{n+1}]/J$ be the quotient
ring, which we note is graded.  Then the residue maps $R^{a(q)}$ to
$F^q/F^{q+1}$. By a theorem of Griffiths \cite{GriffPerRat}, this map is an isomorphism.
For a smooth variety the ``Jacobian ring'' $R$ is finite-dimensional, and so there is
a least integer 
$$
  t = (n+2)(d-2)
  \eqn
  \eqref{topJ}
$$
such that $R^i = 0$ for $i > t$. Moreover, 
and $R^t$ is one-dimensional and the bilinear map
$$
   R^i\times R^{t-i} \map R^t \cong \C.
$$
is a perfect pairing (Macaulay's theorem).  When $R^i$ and $R^{t-i}$
correspond to graded quotients of the Hodge filtration, the pairing corresponds
to the cup product  \cite{CG}.

The derivative of the period map is given by formal differentiation
of the expressions (\xref{ratdiff}).  Thus, if $P_t = P + tQ + \cdots$
represents a family of hypersurfaces and $\omega
= \res(A\Omega/F^\ell)$ represents a family of cohomology classes
on them, then 
$$
  { d \over dt} \res{ A\Omega \over P^{q+1} } = -(q+1) \res { QA\Omega \over P^{q+2}
}.
$$
To show that the derivative of the period map is nonzero,
it suffices to exhibit an $A$ and a $Q$ which are nonzero
in $R$ and such that the product $QA$ is also nonzero.
Here we implicitly use
the identification
$
  T \cong R^d
$
of tangent vectors to the moduli space with the component
of the Jacobian ring in degree $d$.  Thus the natural components
of the differential of the period map,
$$
   T \map \Hom(H^{p,q}(X),H^{p-1,q+1}(X)),
$$
can be identified with the multiplication homomorphism
$$
  R^d \map \Hom(R^a,R^{a+d}),
$$
where $a$ is the degree of the numerator polynomial used in
the residues of the forms (\xref{ratdiff}).  All of these results,
discovered first by Griffiths in the case of hypersurfaces,
hold for weighted hypersurfaces by the results described in 
\cite{Dolgachev} and \cite{Tu}.

Consider now a double cover $Y$ of a hypersurface $X$ of even degree $d$.
If $X$ is defined by $P(x_0 \commadots x_n) = 0$ then $Y$ is defined by
$y^2 +P(x_0 \commadots x_n) = 0$, where $y$ has weight $d/2$ and where the $x$'s
have weight one.  This last equation is homogeneous of degree $d$ with respect
to the given weighting, and $\Omega$ has weight $d/2 + n + 2$.  Thus
$\nu(A,y^2 + P,q)$ is of weight zero if
$a(q) = (q + 1/2)d - (n+2)$.  Since $y$ is in the Jacobian ideal,
we may choose $A$ to be a polynomial in the $x$'s, and we may consider
it modulo the Jacobian ideal of $P$.  Thus the classical considerations
of the residue calculus apply.  If we choose $a(q)$ maximal subject to 
the constraints $p > q$ and $a \ge 0$ then 
$$
 q = \left\{ { n + 1 \over 2 } \right\},
$$
where $\set{ x }$ is the greatest integer {\sl strictly less}
than $x$.  Both conditions are satisfied for $d \ge 4$ except
that for $n = 1$ we require $d \ge 6$.  Thus we have excluded the case $(d,n)
= (4,1)$ in which the resulting  double cover is rational and the period map is
constant. 

Now let $A$ be a polynomial of degree $a$ which is
nonzero modulo the Jacobian ideal.  We must exhibit a polynomial $Q$ of degree
$d$ such that $AQ$ nonzero modulo $J$.  By Macaulay's theorem there is a polynomial
$B$ such that $AB$ is congruent to a generator of $R^t$, hence satisfies
$AB \not\equiv 0 \hbox{ mod $J$}$.  Write $B$ as a linear combination of monomials
$B_i$ and observe that there is an $i$ such that $AB_i \not\equiv 0$.  If $B_i$
is of degree at least $d$, we can factor it as $QB_i'$ with $Q$ of degree $d$. 
Since $AQB_i' \not\equiv 0$, $AQ \not\equiv 0$, as required.

The condition that $B$ have degree at least $d$ reads $a + d \le t$.  Using
the formulas \eqrefer{topJ} for $t$ and the optimal choice for $a$, we see that this
inequality is  satisfied for the range of $d$ and $n$ considered.  This computation
completes the proof of the main theorem in the case $d$ even, $d \ge 4$,
except for the case $(d,n) = (4,1)$.

\section{Rational differentials for higher cyclic covers }
\secref{cycliccoversection}

To complete the proof of the main theorem we must consider 
arbitrary cyclic covers of $\P^{n+1}$ branched along a smooth hypersurface
of degree $d$.  Since the fundamental group of the complement of $X$
is cyclic of order $d$, the number of sheets $k$ must be a divisor of
$d$.  As mentioned in the outline of the proof, there is an automorphism
$\sigma$ of order $k$ which operates on the universal family $\bY$ of 
such covers.  Consequently the local system $\H$ of vanishing cohomology
(cf. \S \xref{zardensitysection})
splits over $\C$ into eigensystems $\H(\mu)$, where $\mu \ne 1$
is a $k$-th root of unity.  Therefore the monodromy representation,
which we now denote by $\rho$, splits as a sum of representations $\rho_\mu$
with values in the groups $\widetilde G(\mu)$ introduced in \eqrefer{tildeGdecomp}.
As noted there we can view $\rho_\mu$ as taking values in a group of linear
automorphisms of $H(\mu)$.  This group is unitary for the
hermitian form $h(x,y) = i^{n+1}(x,\bar y)$ if $\mu$ is
non-real, and that is the case that we will consider here.

Although the decomposition of $\H$ is over the complex numbers, important
Hodge-theoretic data survive.  The hermitian form $h(x,y)$ is
nondegenerate and there is an induced Hodge decomposition, although
$h^{p,q}(\mu) = h^{q,p}(\mu)$ may not hold.  However, Griffiths'
infinitesimal period relation,
$$
   { d \over dt } F^p(\mu) \subset F^{p-1}(\mu)
$$
remains true.  Thus each $\H(\mu)$ is a {\sl complex variation of Hodge
structure}, c.f. \cite{DeMo}, \cite{SimpsonHiggs}.  The
associated period domains are homogeneous for the groups 
$\widetilde G(\mu)$.  

To extend the arguments given above to the unitary representations
$\rho_\mu$ we must extend Deligne's density theorem to this case.
The essential point is that the monodromy groups $\Gamma(\mu)$ are
generated not by Picard-Lefschetz transformations, but by their unitary
analogue, which is a {\sl complex reflection} \cite{Pham},
\cite{Givental}, \cite{Mostow}.   These are linear maps of the form
$$
  T(x) = x \pm (\lambda - 1)h(x,\delta)\delta,
$$
where $h$ is the hermitian inner product defined above,  $h(\delta,\delta) = \pm
1$, where
$\pm$ is the same sign as that of $h(\delta ,\delta )$, and where
$\lambda \ne 1$ is a root of unity.  The vector $\delta$ is an eigenvector of $T$
with eigenvalue $\lambda$ and $T$ acts by the identity on the hyperplane
perpendicular to $\delta$.  It turns out that the eigenvalue $\lambda$ of $T$
is, up to a fixed sign that depends only on the dimension of $Y$, equal
to the eigenvalue $\mu$ of $\sigma$.

In section \xref{unitarydensitysection} we will prove an analogue of Deligne's
theorem \eqrefer{delignedensity} for groups of complex reflections.
It gives the usual dichotomy: either the monodromy group is finite, or it is
Zariski-dense.  In section \xref{complexreflectionsection} we will show that the
monodromy groups $\Gamma(\mu)$ are indeed generated by complex reflections.
It remains to show that the derivative of the period map for the complex
variations of Hodge structure $\H(\mu)$ are nonzero given appropriate
conditions on $d$, $k$, $n$, and $\mu$.

For the computation fix $\zeta = e^{2\pi i/k}$ as a primitive $k$-th 
{\sl root of unity} and let the cyclic action on the universal family
\eqrefer{universalcyclic} be given by $y\circ \sigma = \zeta y$.  
Then the ``volume form'' $\Omega(x,y)$ is an eigenvector
with eigenvalue $\zeta$ and the rational differential
$$
  { y^{i-1} A(x) \Omega(x,y) \over ( y^k + P(x) )^{q+1} }
  \eqn
  \eqref{iratdiff}
$$
has  eigenvalue $\mu = \zeta^i$, as does its residue. 
Thus we will sometimes write $\H(i)$ for 
$\H(\zeta^i)$ and will use the corresponding notations 
$\widetilde G(i)$, $\tilde\rho_i$, etc.  Residues with numerator
$y^{i-1}A(x)$ and denominator $(y^k + P(x))^{q+1}$ span the spaces $H^{p,q}_0(i)$, where
$i$ ranges from 1 to $k-1$.  Moreover, the corresponding space of numerator
polynomials, taken modulo the Jacobian ideal of $P$, is isomorphic via the residue
map to $H^{p,q}_0(i)$.  Since $P$ varies by addition of a polynomial in the $x$'s,
the standard unweighted theory applies to computation of the derivative map.

Let us illustrate the relevant techniques by computing the Hodge numbers and
period map for triple covers of $\P^3$ branched along a smooth cubic surface.
(This period map is studied in more detail in \cite{ACT}.) 
A triple cover of the kind considered is a cubic hypersurface in $\P^4$, and
the usual computations with rational differentials show that $h^{3,0} = 0$, 
$h^{2,1} = 5$.  The eigenspace
$H^{2,1}(i)$ is spanned by residues of differentials with numerator
$ A(x)\Omega(x,y) $ and denominator $(y^3 + P(x))^2$.  Since the degree of
$\Omega(x_0,x_1,x_2,x_3, y)$ is 5, $A$ is must be linear in the variables $x_i$.  Thus
$h^{2,1}(1) = 4$.  The space
$H^{1,2}(1)$ is spanned by residues of differentials with numerator
$ A(x)\Omega(x,y) $ and denominator $(y^3 + P(x))^3$.
Thus the numerator is of degree four, but must be viewed modulo
the Jacobian ideal.  For dimension counts it is enough
to consider the Fermat cubic, whose Jacobian ideal is generated
by squares of variables.  The only square-free quartic in four variables
is $x_0x_1x_2x_3$, so $h^{1,2}(1) = 1$.  Similar computations show
that the remaining Hodge numbers for $H^3(1)$ are zero and
yield in addition the numbers for $H^3(2)$.  One can
also argue that $H^3(1)\oplus H^3(2)$ is defined over $\R$, since the
eigenvalues are conjugate.  A Hodge structure defined over $\R$
satisfies $h^{p,q} = h^{q,p}$.  From this one deduces that
$h^{2,1}(2) = 1$, $h^{1,2}(2) = 4$.  Since there is just one conjugate 
pair of eigenvalues of $\sigma$, there is just one component in the 
decomposition \eqrefer{tildeGdecomp}, $\widetilde G = \widetilde G(\zeta)$,
and this group is isomorphic to $U(1,4)$.  Since the coefficients
of the monodromy matrices lie in the ring $\Z[\zeta]$, where 
$\zeta$ is a primitive cube root of unity, the representation $\tilde\rho$
takes values in a discrete subgroup of $\widetilde G$.  Therefore the complex variation
of Hodge structures define period mappings
$$
   p : U_{3,2} \map B_4/\Gamma',
$$
where $B_4$ is the unit ball in complex 4-space and $\Gamma'$ a
discrete group acting on it.

To show that the period
map $p_i$ is nonconstant it suffices to show that its differential is
nonzero at a single point.  We do this for the Fermat variety. A basis
for $H^{2,1}(1)$ is given by the linear forms $x_i$, 
and a basis for $H^{1,2}$ is given by their product
$x_0x_1x_2x_3$.  Let $m_i$ be the product of all the $x_k$ 
except $x_i$.  These forms constitute a basis for the tangent
space to moduli.  Since $m_ix_i = x_0x_1x_2x_3$, multiplication
by $m_i$ defines a nonzero homomorphism from $H^{2,1}(1)$
to $H^{1,2}(1)$.  Thus the differential of the period map
is nonzero at the Fermat.  In fact it is of rank four, since
the homomorphisms defined by the $m_i$ are linearly independent.  Similar
considerations show that the period map for $\H(2)$ is of rank
four.  The relevant bases are $\set{ y }$ for $H^{2,1}(2)$
 and $\set{ ym_0, ym_1, ym_2, ym_3 }$ for $H^{1,2}(2)$. 
 
 For the general case it will be enough to establish the following.

\proclaim{Proposition.} Let $\bY$ be the universal family of $d$-sheeted
 covers of $\P^{n+1}$ branched over smooth hypersurfaces of degree $d$.
 The derivative of the period map for $\H^{n+1}(1)$ is nontrivial 
 if $n \ge 2$ and $d \ge 3$ or if $n = 1$ and $d \ge 4$.
\endproclaim
\procref{PropositionDerivNonTrivialNgeTwo}

\proof  Elements of
$H^{p,q}(1)$ with $p+q = n+1$ are given by rational differential
forms with numerator $ A(x) \Omega(x,y) $ and denominator $( y^d + P(x) )^{q+1}$.
The numerator must have degree $ a = (q + 1)d - (n+3) $.
As before choose $q$ so that $a$ is maximized subject to the constraints
$p > q$ and $a \ge 0$.  Then $q = \set{ {n / 2} + { 1 / d } }$.
If $n \ge 2$ and $d \ge 3$ or if $n = 1$ and $d \ge 4$, then $a \ge 0$. 
Thus numerator polynomials $A(x)$ which are nonzero modulo the Jacobian ideal
exist.  One establishes the existence of a polynomial $Q(x)$ of degree $d$
such that $QA$ is nonzero modulo the Jacobian ideal using the same argument
as in the case of double covers.

A different component of the period map is required if the branch locus
is a finite set of points, which is the case for the braid group of $\P^1$:

\proclaim{Proposition.} For $n=0$ the period map for $\H^1(i)$ is
non-constant if $d \ge 4$ and $i \ge 2$.
\endproclaim

\proof An element of $H^{1,0}(i)$ is the residue of a rational differential
with numerator $ y^{i-1}A(x_0,x_1)\Omega $ and denominator $ y^d + P(x_0,x_1) $.
The degree of $A$ is $a = d - 2 - i$.  The top degree for the Jacobian
ideal is $2d-4$.  Thus we require $a + d \le 2d - 4$, which is satisfied
if $i \ge 2$.  Since $a \ge 0$, one must also require $d \ge 4$.

We observe that the local systems which occur
as constituents for $k$-sheeted covers, where $k$ divides $d$,
also occur as constituents of $d$-sheeted covers.

\proclaim{Remark.} Let $\H({k,\mu})$ be the complex variation
of Hodge structure associated to a $k$-sheeted cyclic cover
of $\P^{n+1}$ branched along a hypersurface of degree $d$,
belonging to the eigenvalue $\mu$, where $k$ is a divisor
of $d$.  Then $\H({k,\mu})$ is isomorphic to $\H({d,\mu})$.
\endproclaim

\proof Consider the substitution $y = z^{d/k}$ which effects the transformation
$$
  { y^i A(x) \Omega(x,y) \over (y^k + P(x)  )^{q+1} }
 \mapsto
  { (d / k) }{ z^{(i+1)(d/k) -1} A(x) \Omega(x,z) \over 
  ( z^d + P(x) )^{q+1} } .
$$ 
These differentials are eigenvectors with the same eigenvalue.  The
map which sends residues of the first kind of rational differential
to residues of the second defines the required isomorphism.

\section{Complex Reflections}
\secref{complexreflectionsection} 

We now review some known facts on how complex reflections arise for
degenerations of cyclic covers.  When the branch locus acquires
a node, the local equation is 
$$
  y^k + x_1^2 + \cdots + x_{n+1}^2 = t,
\eqn
\eqref{kdoublept}
$$
which is a special case of the situation studied by Pham in \cite{Pham},
where the left-hand side is a sum of powers.  Our
discussion is based on  Chapter 9 of \cite{Milnor} and  Chapter 2 of \cite{Arnold}.

Consider first the case $y^k = t$. It is a family of
zero-dimensional varieties  $\set{ \xi_1(t) \commadots \xi_k(t)}$ whose vanishing
cycles are successive differences of roots,
$$
  \xi_1 - \xi_2, \  
  \ldots,\ 
  \xi_{k-1} - \xi_k,
  \eqn
  \eqref{stdvanishingbasis}
$$
and whose monodromy is given by cyclically shifting indices to the right:
$$
  T( \xi_i - \xi_{i+1} ) = \xi_{i+1} - \xi_{i+2},
$$
where $i$ is taken modulo $k$.  Thus $T$ acts on the $(k-1)$-dimensional 
space of vanishing cycles as a transformation of order $k$.  Over the complex
numbers it is diagonalizable, and the eigenvalues are the $k$-th roots of unity
$\mu \ne 1$.  Note that $T = \sigma_0$ where $\sigma_0$ is the generator for the
automorphism group of the cyclic cover $y^k = t$ given by $y\map \zeta y$, where
$\zeta = e^{2\pi i/k}$ is our chosen primitive $k$-th root of unity.

The intersection product $B$ defines a possibly degenerate bilinear
form on the space of vanishing cycles.  For the singularity $y^k = t$
it is $(\xi_i,\xi_j) = \delta_{ij}$, so relative
to the basis \eqrefer{stdvanishingbasis} it is the negative of the matrix for
the Dynkin diagram $A_{k-1}$ --- the positive-definite matrix with two's along
the diagonal, one's immediately above and below the diagonal, and zeroes elsewhere.

Now suppose that $f(x) = t$ and $g(y) = t$ are families which
acquire an isolated singularity at $t = 0$.  Then $f(x) + g(y) = t$
is a family of the same kind; we denote it by $f \oplus g$.  The theorem of
Sebastiani and Thom \cite{ST}, or 
\cite{Arnold}, cf. Theorem 2.1.3, asserts that vanishing cycles for the
sum of two singularities are given as the join of vanishing cycles for $f$ and
$g$.  Thus, if $a$ and $b$ are vanishing cycles of dimensions $m$ and $n$, then
the join
$a*b$ is a vanishing cycle of dimension $m+n+1$, and, moreover,
the monodromy acts by 
$
  T(a * b) = T(a)*T(b).
$
{}From an algebraic standpoint the join is a tensor product, so one can write
$V(f\oplus g) = V(f)\otimes V(g)$ where $V(f)$ is the space of vanishing
cycles for $f$, and one can write the monodromy operator as
$
  T_{f \oplus g} = T_f\otimes T_g.
$

The {\sl suspension} of a singularity $f(x) = t$ is  by definition the
singularity $y^2 + f(x) = t$ obtained by adding a single square.  If $a$ is
a vanishing cycle for $f$ then $(y_0 - y_1)\otimes a$ is a vanishing
cycle for the suspension, and the suspended monodromy is given by
$$
  T( (y_0 - y_1)\otimes a ) = - (y_0 - y_1)\otimes T(a) .
$$
In particular, the local monodromy of a singularity
and its double suspension are isomorphic.

The intersection matrix $B'$ of
a suspended singularity (relative to the same canonical basis) is a function of
the intersection matrix $B$ for the given singularity, cf. Theorem 2.14 of
\cite{Arnold}.  When the bilinear form for $B$ is symmetric, the rule for producing $B'$ from $B$
is: make the diagonal entries zero and change the sign of the above-diagonal entries.  When $B'$
has an even number of rows of columns, the determinant is one, and when the number
of rows and columns is odd, it is zero.  Thus the intersection matrix
for $x^2 + y^k = t$ is nondegenerate if and only if $k$ is odd.  In addition,
the intersection matrix of a double suspension is the negative of the given
matrix.  Thus the matrix of any suspension of
$y^k = t$ is determined.  It is nondegenerate if the  dimension of the
cyclic cover \eqrefer{kdoublept} is even or if the dimension is odd and $k$ is
also odd.  Otherwise it is degenerate.

It follows from our discussion that the space of vanishing cycles
$V$ for the singularity \eqrefer{kdoublept} is $(k-1)$-dimensional
and that the local monodromy transformation is
$
  T = \sigma_0\otimes(-1)\otimes\cdots\otimes(-1)
$
where $\sigma_0$ is the covering automorphism $y\map \zeta y$ for $y^k = t$.  Thus
$T$ is a cyclic transformation of order $k$ or $2k$, depending on whether
the dimension of the cyclic cover is even or odd.  In any case, $T$ is
diagonalizable with eigenvectors $\eta_i$ and eigenvalues $\lambda_i$,
where $\lambda_i = \pm \mu_i$ with $\mu_i = \zeta^i$ where $\zeta$ is our fixed 
primitive 
$k$-th root of unity and $i = 1,\cdots , k-1$.
  Note that the cyclic automorphism $\sigma$ of the universal family
\eqrefer{universalcyclic}, given by 
$y\mapsto\zeta y$  
acts as  $\sigma_0\otimes(+1)\otimes\cdots\otimes(+1)$ on the vanishing homology
of \eqrefer{kdoublept}. Thus the eigenspaces of $\sigma$ and $T$ coincide, and
their respective eigenvalues differ by the fixed sign $(-1)^{n+1}$.  Since the
eigenvalues
$\mu_i$ are distinct, the eigenvectors $\eta_i$ are orthogonal with respect
to the hermitian form.  Thus $h(\eta_i,\eta_i) \ne 0$.  Moreover the sign of 
$h(\eta_i,\eta_i)$ depends only on the index $i$, globally determined on
\eqrefer{universalcyclic}, independently of the particular smooth point on 
the discriminant locus whose choice is implicit in \eqrefer{kdoublept}.   We
conclude that on the space of vanishing cycles,
$$
  T(x) = \sum_{i = 1}^{k-1} \lambda_i{h(x,\eta_i) \over
h(\eta_i,\eta_i)}\eta_i,
  \eqn
  \eqref{TVcxreflectionformula}
$$
where $\lambda_i = (-1)^{n+1} \mu_i$.

Now consider a cycle $x$ in $H^{n+1}(Y_{\tilde o})$, and {\sl suppose that}
$k$ {\sl is odd}.  Then the intersection form on the space $V$ of local vanishing
cycles
for the degeneration \eqrefer{kdoublept} is {\sl nondegenerate}. Consequently
$H^{n+1}(Y_{\tilde o})$ splits orthogonally as
$V
\oplus V\perp$.  The action on $H^{n+1}(Y_{\tilde o})$ of the monodromy
transformation $T$ for the meridian corresponding to the degeneration
\eqrefer{kdoublept}   
 is given by \eqrefer{TVcxreflectionformula} on $V$ and 
by the identity on $V\perp$.  Thus it is given for arbitrary $x$ by the formula
$$
  T(x) = x + \sum_{i=1}^{k-1} (\lambda_i-1){h(x,\eta_i) \over
h(\eta_i,\eta_i)}\eta_i .
  \eqn
  \eqref{Tcxreflectionformula}
$$
Finally, for each $i = 1,\cdots ,k-1$ we can normalize the eigenvector $\eta_i$
to an eigenvector $\delta_i$ satisfying $h(\delta_i , \delta_i ) = \epsilon_i = 
\pm 1$.  To summarize, we have proved the following:

\proclaim{Proposition.} Consider the family \eqrefer{universalcyclic} of $k$-fold
cyclic covers  of $\P^{n+1}$ branched over a smooth hypersurface of degree $d$,
where both $k$ and $d$ are odd.  Let $T$ be the monodromy corresponding to
a generic degeneration of the branch locus, as in \eqrefer{kdoublept}.
Then $T$ acts on the $i$-th eigenspace of the cyclic automorphism  
$\sigma$ (defined by $y\mapsto\zeta y$ in \eqrefer{universalcyclic}) by a
complex reflection with eigenvalue
$\lambda_i = (-1)^{n+1} \zeta^i$.  Thus 
$$
T(x) = x + \epsilon_i (\lambda_i - 1) h(x,\delta_i)\delta_i
$$
holds for all $x\in\H (i)$. 
\endproclaim

\proclaim{Remark.} \rm In remark \xref{remarkMeridiansInfiniteOrder}.a we observed that
the meridians of $\Phi_{d,n}$ are of infinite order for $n$ odd and for 
$n$ even, $d \ge 4$ even.  Consider now the case in which $n$ is even and $d$ is odd,
let $\zeta = \exp(2\pi i/d)$, and let $\bar\rho'$ be the
corresponding representation, in which meridians of $\widetilde \Delta$
correspond to complex reflections of order $2d$. These complex reflections
and their powers different from the identity are non-central if the $\zeta$ eigenspace
has dimension at least two, which is always the case for $d \ge 3$, $n \ge 2$.
Thus $\bar\rho'(\gamma)$ has order $2d$.  By this simple argument
we conclude that in the stated range of $(n,d)$, meridians always have order greater than two. 
However, our argument does not give the stronger  result \xref{remarkMeridiansInfiniteOrder}.b
asserted by Kontsevich.

\endproclaim

\section{Density of unitary monodromy groups}
\secref{unitarydensitysection}

 
We now show how the argument  Deligne used in \cite{DeWeTwo}, section 4.4,
to prove Theorem \xref{delignedensity} can
be adapted to establish a density theorem for groups generated
by complex reflections on a space $\C(p,q)$ endowed with
a hermitian form $h$ of signature $(p,q)$.  If $A$ is a subset of $\C(p,q)$ or 
of $U(p,q)$, we use $PA$ to denote its projection in $\P(\C(p,q))$ or $PU(p,q)$.

\proclaim{Theorem.} Let $\epsilon = \pm 1$ be fixed, and let $\Delta$ be a set of
vectors in a hermitian space
$\C(p,q)$ which lie in the unit quadric $h(\delta,\delta) = \epsilon$.
Fix a root of unity $\lambda \ne \pm 1$ and let $\Gamma$ be the 
subgroup of $U(p,q)$ generated by the complex reflections 
$s_\delta(x) = x + \epsilon ( \lambda - 1 )h(x,\delta)\delta$ for all $\delta$ in
$\Delta$.  Suppose that $p+q >1$, that $\Delta$ consists of a single
$\Gamma$-orbit, and that $\Delta$ spans $\C(p,q)$.  Then either $\Gamma$ is 
finite or $P\Gamma$ Zariski-dense in $PU(p,q)$.
\endproclaim
\procref{udensitytheo}

Let $\bar\Gamma$ be the Zariski closure of a subgroup $\Gamma$
of $U(p,q)$ which contains the $\lambda$-reflections 
for all vectors $\delta$ in a set $\Delta$.  Then 
$\bar\Gamma$ also contains the $\lambda$-reflections for
the set $R = \bar\Gamma\Delta$.  Indeed, if $g$ is an element
of $\bar\Gamma$, then 
$$
  g^{-1}s_\delta g = s_{g^{-1}(\delta)}.
  \eqn\eqref{reflectionconjugacy}
$$
Thus it is enough to establish the following result in order
to prove our density theorem:

\proclaim{Theorem.} Let $\epsilon = \pm 1$ be fixed, and let 
$R$ be a set of vectors in a hermitian space
$\C(p,q)$ which lie in the unit quadric $h(\delta,\delta) = \epsilon$.
Fix a root of unity $\lambda \ne \pm 1$ and let $M$ be the smallest algebraic 
subgroup of $U(p,q)$ which contains the complex reflections 
$s_\delta(x) = x + \epsilon ( \lambda - 1 )h(x,\delta)\delta$ 
for all $\delta$ in $R$. 
Suppose that $p+q >1$, that $R$ consists of a single $M$-orbit,
and that $R$ spans $\C(p,q)$.  Then either $M$ is finite or
$PM = PU(p,q)$.
\endproclaim
\procref{udensityprop}

We begin with a special case of the theorem for groups
generated by a pair of complex reflections.

\proclaim{Lemma.} Let $\lambda \ne \pm 1$ be a root of unity, and let 
 $U$ be the unitary group of a nondegenerate
hermitian form on $\C^2$.  Let $\delta_1$ and $\delta_2$ be
independent vectors with nonzero inner product, and let 
$\Gamma$ be the group generated by complex reflections with common
eigenvalue $\lambda$. Then either $\Gamma$ is finite or its
image in the projective unitary group is Zariski-dense.
In the positive-definite case $\Gamma$ is finite if and only if
the inner products $(\delta_1,\delta_2)$ lie in a 
fixed finite set $S$ which depends only on $\lambda$ and $h$.  In the indefinite
case
$\Gamma$ is never finite.
\endproclaim
\procref{Ulemma}

We treat the definite case first.  To begin, note
that the group $U$ acts on the Riemann sphere $\P^1$ via the natural map
$U \map PU$, where $PU$ is the projectivized unitary group.  Let $PR$ be
the image of $R \subset \C^2$ in $\P^1$.  Since $\lambda$ is a root of 
unity, the projection $P\Gamma$ is a
finite group if and only if $\Gamma$ is.  The finite subgroups of rotations
of the sphere are well known.  There are two infinite series: the cyclic
groups, where the vectors $\delta$ are all proportional, and the dihedral
groups where $\lambda = -1$.  There are three additional groups, given
by the symmetries of the five platonic solids, and  $S$ is the set of possible
values of $h(\delta_1,\delta_2)$ that can arise for these three groups.

We suppose that
$(\delta_1,\delta_2)$ lies outside $S$, so that $P\Gamma$ is infinite.
Then its Zariski closure $PM$ is either $PU$ or a group whose identity
component is a circle.  In this case $PR$ contains a great circle $\alpha$.
However, $PR$ is stable under the action of $PM$, hence under the rotations
corresponding to axes in $PR$.  Since $\lambda \ne \pm 1$, the orbit
$PR$ contains additional great circles which meet $\alpha$ in an angle $0 < \phi
\le \pi/2$.  The union of these, one for each point of the given
circle, forms a band about the equator, hence has nonempty interior. Such a
set is Zariski-dense in the Riemann sphere viewed as a real algebraic
variety. Since
$PR$ is a closed real algebraic set, $PR = S^2$.  Since
$PR \cong PM/H$, where $H$ is the isotropy group of a point on the sphere,
$PM = PU$.

In the case of an indefinite hermitian form, the group
$U =  U(1,1)$, acts on the hyperbolic plane via the projection
to $PU$, and $P\Gamma$
is a group generated by a pair of elliptic elements of 
equal order but with distinct fixed points.  One elliptic
element moves the fixed point of the other, and so their commutator $\gamma$
is hyperbolic (c.f. Theorem 
7.39.2 of \cite{Beardon}).  The Zariski closure of the cyclic group
$\set{\gamma^n}$ is a one-parameter subgroup of $PU$. 
Consequently the orbit $PR$ contains a geodesic $\alpha$ through one of the
elliptic fixed points. By \eqrefer{reflectionconjugacy} the other points of
$\alpha$ are
 fixed points of other elliptic transformations in $PM$.  Now the orbit $PR$
contains the image of  $\alpha$ under each of these transformations,
and so $PR$ contains an open set of the hyperbolic plane.  This implies that 
either $PM = PU$ or $PM$ is contained in a parabolic subgroup.  Since $PM$
contains non-trivial elliptic elements that last possibility cannot occur,  
and so  $PM = PU$.

Next we show that if the set $R$ which defines the reflections
is large, then so is the group containing those reflections.

\proclaim{Lemma.} Fix a root of unity $\lambda\ne\pm 1$ and $\epsilon = \pm 1$.  
Let $R$ be a semi-algebraic subset of the unit quadric 
$h(\delta,\delta) = \epsilon$.
Let $M$ be the smallest algebraic subgroup of $U(p,q)$
containing the complex reflections  $s_\delta(x) = x + \epsilon (\lambda - 1)
h(x,\delta)\delta$, $\delta\in R$.  If $p+q >1$ and if $PR$ 
is Zariski-dense in $\P(\C(p,q))$,
then
$M =PU(p,q)$.
\endproclaim

The proof is by induction on $n = p+q$.  For $n = 2$ the result
follows from the proof of lemma \xref{Ulemma}.  Let $n>2$ and assume 
$p\le q$. Then $q\ge 2$.  Fix a codimension two subspace of $\C (p,q)$ of 
signature $(p,q-2)$ and let $W_t$ be the pencil of hyperplanes of $\C (p,q)$ 
containing this codimension two subspace.  Then the restriction of $h$ to each 
$W_t$ is a non-degenerate form of signature $(p,q-1)$. 

Consider a subgroup $M$ of $U(p,q)$ which satisfies
the hypotheses of the lemma, and let $R_t = R\cap W_t$.  Since $PR$, respectively 
$PR_t$ is semi-algebraic in $\P (\C (p,q))$, respectively in $PW_t$, it is 
Zariski dense if and only if it has non-empty interior in the analytic topology.  
Thus $R$ has non-empty interior in $\P(\C(p,q))$, and so for dimension reasons 
$PR_t$ has non-empty interior in $PW_t$ for generic $t$.  Thus $PR_t$ is 
Zariski dense in $PW_t$ for generic $t$.

Fix one such value of $t$, let $W = W_t$ and let $M'(R\cap W)\subset M$ denote the 
smallest algebraic subgroup of $M$ containing $R\cap W$.  Let $M(R\cap W)$ denote 
the set of restrictions of elements of $M'(R\cap W)$ to $W$.  Then $R\cap W$ and 
$M(R\cap W)$ satisfy the induction hypothesis, thus $PM(R\cap W) = PU(W)$.
 Now the orthogonal
complement of $W$ is a Zariski closed set, as is $W \cup W\perp$.
Since $R$ is Zariski-dense there is a $\delta$ in $R - W$ 
and a $\delta'$ in $W$ such that $h(\delta,\delta') \ne 0$. 
Consider the function $f_\delta(x) = h(x,\delta')$.  If it is constant
on the Zariski closure $C$ of $R\cap W$, then the derivative
$df_\delta$ vanishes on $C$.  Therefore $C$ lies in the intersection
of the hyperplane $df(x) = 0$ with $W$, which is a proper
algebraic subset of $W$.  Consequently $R \cap W$ is not Zariski-dense,
a contradiction.  Thus $f_\delta$ is nonconstant and so we can choose
$\delta$ in $R \cap W$ such that $h(\delta',\delta)$ lies outside
the fixed set $S$.  Then lemma \xref{Ulemma} implies that the unitary group  of
the plane $F$ spanned by $\delta$ and $\delta'$ is contained in $M$.  But $U(W)$
and $U(F)$ generate $U(p,q)$ and the proof of the lemma is complete.

To complete the proof of Theorem (\xref{udensityprop})
we must show that either $R$ is sufficiently large or that $M$
is finite.  Observe that since $R$ is an $M$-orbit, it is a semi-algebraic set.  
Let $W$ be a subspace of $\C(p,q)$
which is maximal with respect to the property ``$W \cap R$
is Zariski-dense in the unit quadric of $W$.''  Our aim is to show
that either $W = \C(p,q)$ or that $M$ is finite.  Consider 
first the case $W = 0$.  Then the inner products
$h(\delta,\delta')$ for any pair of elements in $R$ lie in the 
fixed finite set $S$ of lemma \xref{Ulemma}. Now let $\delta_1 \commadots
\delta_n$ be a basis of $\C(p,q)$ whose elements are chosen from $R$.
Then the inner products $h(\delta,\delta_i)$ lie in $S$
for all $\delta$ and $i$.  Consequently $R$ is a finite set
and $M$, which is faithfully represented as a group of permutations
on $R$, is finite as well.

Henceforth we assume that $W$ is nonzero.  If it is not maximal
there is a vector $\delta$ in $R - W$ and we may consider
the function  $f_\delta(x) =  h(x,\delta)$ on the set $R \cap W$.  If $f_\delta$
is identically zero for all $\delta$ in $R - W$, then $R \subset W \cup W\perp$.
Therefore $\C(p,q) = W + W\perp$, from which one concludes that
$W = W \oplus W\perp$ and so $M$ is a subgroup of $U(W)\times
U(W\perp)$. But $R$ consists of a single $M$-orbit and contains a point of $W$,
which implies that $R \subset W$, a contradiction.

We can now assume that there is a $\delta \in R - W$ such that
the function $f_\delta$ is not identically  zero.
If one of these functions is not locally
constant, then it must take values outside the set $S$.  Then
the inner product $(x,\delta)$ lies outside $S$ for an open dense
set of $x$ in $R \cap W$.  For each such $x$, $R$ is dense in 
the span of $x$ and $\delta$.  We conclude that $R$ is dense in
$W + \C\delta$.  Thus $W$ is not maximal, a contradiction.

At this point we are reduced to the case in which all the functions
$f_\delta$ are locally constant, with at least one which is not 
identically zero. To say that $f_\delta$ is locally constant
on a dense subset of the unit quadric in $W$ is to say that its
derivative is zero on that quadric.  Equivalently, tangent spaces
to the quadric are contained in the kernel of $df_\delta$,
that is, in the hyperplane $\delta\perp$.  But if all tangent spaces
to the quadric are contained in that hyperplane, then so is the quadric
itself.  Then the function in question is identically zero, contrary
to hypothesis. The proof is now complete.

To apply the density theorem we need to show that the
``complex vanishing cycles'' contain a basis for the vanishing cohomology and 
form
a  single orbit.   These cycles are by
definition the eigencomponents of ordinary vanishing cycles.  
Consider now a generalized Picard-Lefschetz transformation given by
\eqrefer{Tcxreflectionformula}.  It can be rewritten as
$$
   \rho(\gamma)(x) = 
   x + \sum \epsilon_i(\lambda_i-1) h(x,\delta_i) \delta_i,
$$
where the $\delta_i$ are complex vanishing cycles and the $\lambda_i$
are suitable complex numbers.  Let
$$
   \rho(\gamma')(x) = x + \sum\epsilon_i (\lambda_i-1) h(x,\delta'_i) \delta'_i
$$
be another generalized Picard-Lefschetz tranformation.  If
$\gamma' = \kappa^{-1}\gamma\kappa$ then the two preceding equations
yield
$$
  \sum \epsilon_i(\lambda_i - 1)  h(x,\delta'_i) 
\delta'_i 
    = 
  \sum \epsilon_i(\lambda_i - 1 ) h(\kappa.x,\delta_i) 
  \kappa^{-1}.\delta_i ,
$$
where $\kappa.x$ stands for $\rho(\kappa)(x)$.  
Comparing eigencomponents on each side we find
$$
   \delta'_i = \kappa^{-1}.\delta_i,
$$
as required.    By the same argument as used in 
\S \xref{zardensitysection}, one sees that the complex vanishing cycles
span $H(i)$.

\section{Bounds on the real and complex rank}
\secref{rankbounds}

In this section we derive lower bounds for the complex and real ranks of the
groups  
$G_{d,n}$  of  automorphisms of the primitive cohomology
$H^n_o(X_{d,n},\R)$ where
$X_{d,n}$ is a hypersurface of degree $d$ and dimension $n$.  
 Recall that for a field $k$,
the $k$-rank is the dimension of the largest subgroup that
can be diagonalized over $k$.  These bounds
complete the outline of proof.  We also show that all the eigenspaces 
of the cyclic automorphism $\sigma$ have the same dimension.  

The main result is the following:

\proclaim{Lemma}.  The complex rank of $G_{d,n}$ is at least 
five
for $d \ge 3$, $n \ge 1$, with the exception of $(d,n) = (3,1)$, for which it is
 one, and $(d,n) = (4,1), (3,2)$ for which it is three. Under the same
conditions the real rank is at least two with the exception of the cases $(d,n)
= (3,1),\ (3,2)$ for which the real ranks are one and zero, respectively.
\endproclaim
\procref{rankboundslemma}

To prove the first assertion we note that the complex rank is given by 
$\rank_\C G_{d,n} = [ B_{d,n} / 2 ]$ 
where $[x]$ is the greatest integer in $x$ and where 
$ B_{d,n} = \dim H^n_o(X_{d,n})$ is the
primitive middle Betti number. To  compute it we compute the Euler
characteristic $\chi_{d,n}$ recursively using the fact that a $d$-fold cyclic
cover of $\P^{n}$ branched along a hypersurface of degree $d$ is a hypersurface
of degree $d$ in $\P^{n+1}$. Thus, mimicking the  proof of Hurwitz's formula for
Riemann surfaces, we have
$$
   \chi_{d,n} = d\,\chi(\P^n - B) + \chi(B) = d(n+1) + (1-d)\chi_{d,n-1} .
$$
Since $\chi_{d,0} = d$, the Euler characteristics of all hypersurfaces are 
determined.  Rewriting this recursion relation in terms of the $n$-th
primitive Betti number  we obtain 
$$
   B_{d,n} = (d-1)\left(B_{d,n-1} + (-1)^n\right),
   \eqn
   \eqref{bettinumberrecursion}
$$
{}From it we deduce an expression in closed form:
$$
    B_{d,n} = (d-1)^n \,(d-2) + { ( d-1 )^n - (-1)^n \over d  } + (-1)^n .
    \eqn
   \eqref{bettinumberformul}
$$
{}The preceding two formulas
imply that $B_{d,n}$ is an increasing
function of $n$ and of $d$.  Now assume $d \ge 3$, $n \ge 1$. Then $d+n\ge 4$.
If $d+n \le 6$, then $(d,n) = (3,3), (4,2), (5,1)$ and $B_{3,3} = 10 , B_{4,2} =
21, B_{5,1} = 12$. Thus $B_{d,n} \ge 10$ except when $d+n = 4$
or $5$.  These are the cases $(d,n) = (3,1), (3,2), (4,1)$ where $B_{d,n} =
2,6,6$ respectively.  The inequalities on the complex rank are now established.

Let us now turn to the proof of the second assertion of the lemma. For $n$ odd
the group $G_{d,n}$ is a real symplectic group.  Its real and complex ranks
are the same, and so the bound follows from the first assertion.  For $n$ even
the group
$G_{d,n}$ is the orthogonal group of the cup product on the primitive
cohomology.  This bilinear form has signature
$(r,s)$, and the real rank of $G$ is the minimum of $r$ and $s$. The signature is
computed from the Hodge decomposition: 
$r$, the number of positive eigenvalues, is the sum of the 
$h^{p,q}$ for $p$ even, while $s$ is the sum for $p$ odd.   According to the
first inequality of lemma \xref{hodgeinequalities}, the Hodge numbers
$h^{p,q}(d,n)$ of $X_{d,n}$ satisfy 
$h^{p,q}(d+1,n) > h^{p,q}(d,n)$.  Thus the real rank is an increasing function of
the degree.  Consequently it is enough to show that it is at least two for
quartic surfaces and for cubic hypersurfaces of dimension four or more.  For
quartic hypersurfaces $h^{2,0} = 1$ and $h^{1,1} = 19$, so $(r,s) = (2,19)$. 
For cubic hypersurfaces there is a greatest integer
$p \le n$ such that $h^{p,q} \ne 0$, where $p+q = n$.   We will compute this
``first'' Hodge number and see that under the hypotheses of the lemma, $p > q$. 
Since $n$ is even,
$h^{p,q}$ and $h^{q,p}$ have the same parity.  Thus one of $r$, $s$ is at least
two.  According to the second inequality of lemma 
\xref{hodgeinequalities}, $h^{p-1,q+1}(d,n) > h^{p,q}(d,n)$ if $p > q$. Thus
$h^{p-1,q+1}(d,n) > h^{p,q}(d,n) > 0$.  We conclude that the other component of
the signature, $s$ or $r$, must be at least two.  For the Hodge numbers of cubic
hypersurfaces of dimension $n = 3k + r$ where $r = 0$, 1, or 2, one uses the
calculus of \cite{GriffPerRat} to show the following:  
(a) if $n \equiv 0 \mod 3$ then the first
Hodge number is $h^{2k,k} = n+2$, 
(b) if $n \equiv 1 \mod 3$ then it is $h^{2k+1,k} = 1$,
(c) if $n \equiv 2 \mod 3$ then it is $h^{2k+1,k+1} = (n+1)(n+2)/2$.
When $k > 0$ these Hodge numbers satisfy $p > q$, and so 
the proof of the lemma is complete.

\proclaim{Lemma}
 Let $h^{p,q}(d,n)$ be the dimension of $H^{p,q}_o(X_{d,n})$. Then
the inequalities below hold:
$$
\eqalign{
  & h^{p,q}(d+1,n) > h^{p,q}(d,n) \cr
  & h^{p,q}(d,n) > h^{p+1,q-1}(d,n) \hbox{ if $p \ge q$} \cr
}
$$
\endproclaim
\procref{hodgeinequalities}

\proof{}  It is enough to prove the inequalities when $X_{d,n}$
is the Fermat hypersurface defined by 
$F_d(x) = x_0^d  + \cdots + x_{n+1}^d = 0$.  Because of the 
symmetry $h^{p,q} = h^{q,p}$, it is also enough to prove
the inequalities for $p \ge q$.  To this end recall that
$h^{p,q} = \dim R^a$,
where $R$ is the Jacobian ring for $F_d$ and where $a = (q+1)d - (n+2)$
is the degree of the adjoint polynomial in the numerator
of the expression
$$
  \res { A \Omega \over F_d^{q+1} }.
$$
Now there is a map $\mu: R^{a(q,d)}(F_d) \map R^{a(q,d+1)}(F_{d+1})$ defined 
by $\mu(P) =  (x_0 \cdots x_q) P$.  This makes sense because
$q \le n$.  We claim that that resulting map from $H^{p,q}(X_{d,n})$ to
$H^{p,q}(X_{d+1,n})$ is injective but not surjective.

To prove the claim, observe that the Jacobian ideal is generated by 
the powers $x_i^{d-1}$ and so
has a vector space basis consisting of monomials $x^M$. The same
is true of the quotient ring $R(F_d)$.  
Indeed, a basis is given by (the classes of) 
those monomials not divisible by $x_i^{d-1}$
for any $i$.  Now consider a polynomial which represents
an element of the kernel of $\mu$.
It can be be chosen to be a linear combination of monomials $x^M$ which are not
divisible by $x_i^{d-1}$ for any $i$.  Its image is represented by a linear
combination of monomials $(x_0 \cdots x_q)x^M$.   Each of these is divisible by some
$x_i^d$.  Thus either $x^M$ is divisible by $x_i^d$, $i > q$, a contradiction,
or by $x_i^{d-1}$, $i \le q$, also a contradiction.  Thus injectivity part the claim
is  established. 

For the surjectivity part note that image of the map $\mu$ has a basis of monomials
$x^M$ which are divisible by $x_i$ for $i = 0 \commadots q$.  Thus, to show that
$\mu$ is not surjective it suffices to show that there is a monomial for
$R^{a(q,d+1)}(F_{d+1})$  that is not divisible by $x_0$.  Such a monomial has the
form
$x_1^{M_1} \cdots x_{n+2}^{M_{n+2}}$ where $M_i \le d-1$.  It 
exists if $a(q,d+1) \le (n+1)(d-1)$.  The largest relevant values of 
$q$ and $a(q,d+1)$ are $n/2$ and $(n/2 + 1)d - (n+2)$.  For these the
preceding inequality holds and so the first inequality
of the lemma holds strictly.

For the second inequality we use the fact that basis elements for the 
Jacobian ring of $F_d$ correspond to lattice points of the cube in $(n+2)$-space 
defined by the inequalities $0 \le m_i \le d-2$.  A basis for $R^a$
corresponds to the  set of lattice points which lie on the convex subset $C(a)$
of the cube obtained by slicing it with the hyperplane $m_0 + \cdots + m_{n+1} =
a$. The volume of $C(a)$ is a strictly increasing function of $a$ for
$0 \le a \le t/2$, where $t = (n+2)(d-2)$.  For $t/2 \le a \le t$ the volume
function $V(a)$ is strictly decreasing, and in general its graph is symmetric around
$a = t/2$. Let $L(a)$ be the number of lattice points in $C(a)$. If $L(a)$ satisfies
the same monotonicity properties as does $V(a)$, then the second inequality
follows.  To show this, we prove the following result.

\proclaim{Lemma.} Let $L_{d,n}(k)$ be the number of points in the set
$\LL_{d,n}(k) = \sett{ x \in \Z^n }{ 0 \le x_i \le d,\ x_1 + \cdots + x_n = k }$.
Assume that $n > 1$.  Then $L_{d,n}(k)$ is a strictly increasing function of $k$ for $k < dn/2$ 
and is symmetric around $k = dn/2$.
\endproclaim

\proof{} Symmetry follows from the bijection 
$\LL_{d,n}(k) \map \LL_{d,n}(dn - k)$ given by  $x \mapsto \delta - x$
where $\delta = (d \commadots d)$.  We shall say that these two sets
are dual to eachother.  For the inequality we argue by induction, noting
first that
$L_{d,2}(k) = k + 1$ for $k \le d$.  Now observe that $\LL_{d,n}(k)$
can be written as a disjoint union of sets 
$S_i = \sett{ x \in \LL_{d,n}(i) }{ x_n = k-i }$ where $i$ ranges from 
$k-d$ to $k$.  Thus
$$
  L_{d,n}(k) = \sum_{i = k-d}^k L_{d,n-1}(i) .
$$
Consequently
$$
  L_{d,n}(k+1) - L_{d,n}(k) =  L_{d,n-1}(k+1) - L_{d,n-1}(k-d) .
$$
By the induction hypothesis the right-hand side is positive if
$k-d < (n-1)d/2$ and if $k+1$ is not greater than the index dual to $k-d$,
namely $(n-1)d - (k-d)$.  Thus we require also that $k+1 \le (n-1)d - (k-d)$.
Both inequalities hold if $k < nd/2$, which is what we assume.  Thus the proof
is complete.


\subheading{Dimension of the eigenspaces.}

We close this section by noting that the eigenspaces $H^n(X)(\lambda)$
for $\lambda\ne 1$ all have the same dimension, explaining why the
primitive middle Betti number is divisible by $d-1$, where $d$ is the
degree.  Indeed, we have the following,
$$
  \dim H^n(X,\C)(\lambda) 
    = \dim H^n(X,\C)(\mu)
      = \dim H^n(\P^n - B,\C) + (-1)^n .
   \eqn
   \eqref{eigenspacedimform}
$$
When the degree is prime there is a short proof: consider the field $k
= \Q[\omega]$ where $\omega$ is a primitive $d$-th root of unity and
observe that its Galois group permutes the factors $H^n(X,k)(\lambda)$
for $\lambda \ne 1$.  For the general case let $p: X \map \P^n$ be the
projection and note that $H^n(X,\C) = H^n(\P^n, p_*\C)$.  The group of
$d$-th roots of unity acts on $p_*\C$ and decomposes it into
eigensheaves $\C_\lambda$, where $\lambda^d = 1$. Thus the
$\lambda$-th eigenspace of $H^n(X,\C)$ can be identified with
$H^n(\P^n,\C_\lambda)$.  The component for $\lambda = 1$ is
one-dimensional and is spanned by the hyperplane class.  For $\lambda
\ne 1$ the sheaf $\C_\lambda$ is isomorphic to the extension by zero
of its restriction to $\P^n - B$.  Thus the eigenspace can be
identified with $H^n(\P^n - B, \C_\lambda)$. By the argument of
lecture 8 in
\cite{CKM} used in the proof of  vanishing theorems, the groups $H^i(\P^n - B,
\C_\lambda)$ vanish for
$i \ne n$, $\lambda \ne 1$.  Thus $\dim H^n(\P^n - B,\C_\lambda) =
(-1)^n \chi(\lambda)$, where $ \chi(\lambda)$ is the Euler
characteristic of $\C_\lambda$.  Fix a suitable open tubular
neighborhood $U$ of $B$ and a good finite cell decomposition $K$ of
$\P^n - U$.  Then $\chi(\lambda)$ is the Euler characteristic of the
complex of $\C_\lambda$-valued cochains on $K$, which depends only on
the number of cells in each dimension, not on $\lambda$.  This
establishes the first equality above.  For the second use
$\chi(\lambda) = \chi(1)$ and the vanishing of $H^i(\P^n - B,\C)$ for
$i \ne n, 0$.


\section{Remarks and open questions}
\secref{remarkssection}

We close with some remarks on (a) the possiblity of an isomorphism $\Phi \cong
\Gamma\times\Gamma'$, (b) the impossibility of producing additional
representations by iterating the suspension (globally), and (c) generalizations of
the main theorem.

\subheading{(A) Products}

So far everything that has been said is consistent with
an isomorphism between $\Phi$ and the product 
$\Gamma\times\Gamma'$, where $\Gamma'$ is the monodromy group
$\bar \rho'(\Phi)$.  This, however, is not
the case, at least for surfaces, for we can show that {\sl if $k$ is a divisor of $d$ and $d$ is
odd,  then $\Phi_{d,2}$ and $\Gamma\times\Gamma'$ are not isomorphic}.
The argument is based on the fact that the abelianization
of $\Phi$ is a cyclic group of order equal to the
degree of the discriminant, which we denote by $r$.
This is because (a) the generators $g_1 \commadots g_r$ of $\Phi$ are mutually
conjugate, hence equal in the abelianization, (b) $g_1 \cdots g_r = 1$, (c)
the additional relations are trivial when abelianized.  See \cite{Zar}.
For the last point note that $\Phi$ is also the fundamental group
of the complement of a generic plane section $\Delta'$
of $\Delta$.  This complement has nodes and cusps as its
only singularities.  The nodes yield relations of the 
form $gg' = g'g$ where $g$ and $g'$ are conjugates of the
given generators.  The cusps yield braid relations
$gg'g = g'gg'$.  Both are trivial in the abelianization.
Thus the abelianization is generated by a single
element with relation $g^r = 1$.  The degree of 
the discriminant is given in \cite{DolgLib}, page 6, line 2:
$$
  r = \hbox{deg}(\Delta) = 4(d-1)^3.
$$
If $\Phi$ is isomorphic to the Cartesian product,
then there is a corresponding isomorphism of abelianizations.
Let us therefore compute what we can of
the abelianizations of $\Gamma$ and $\Gamma'$.   
For $\Gamma$ we note that the generators are the elements
$g_i$ as above satisfying additional relations which include
$g_i^2 = 1$.  Therefore  $\Gamma$ abelianized is a quotient of $\Z/2$.
Consider next the case of $\Gamma'$ for cyclic covers of 
degree $k$.  Then $\Gamma'$ 
is a product of groups
$\Gamma'(i)$ for $i = 1 \commadots k-1$.  Generators
and relations are as in the previous case except that
among the additional relations are $g_i^{2k} = 1$ instead of $g_i^2 = 1$.
Therefore the abelianization is a quotient of $\Z/{2k}$.
Consequently the abelianization of the
product $\Gamma\times \Gamma'$ is a quotient
of the product of $\Z/2$ with a product of $\Z/2k$'s.  
But the largest the order of an  element in such
a quotient can be is $2k$, which is always 
less than the degree of the discriminant,
provided that $d > 2$, which is the case.

\subheading{(B) Suspensions}

Since $\Phi$ is not in general isomorphic to 
$\Gamma\times\Gamma'$ it is natural to ask
whether there are further representations
with large kernels.  One potential
construction of new representations is
given by iterating the suspension.  By this we mean that
we take repeated double covers.  Unfortunately,
this produces nothing new, since it turns
out that the global suspension is periodic of period two.
To make a precise statement, let $P(x)$ be a polynomial
of degree $2d$ which defines
a smooth hypersurface $X$ in $\P^n$.
Let $X(2)$ be the hypersurface defined by
$$
   P(x) + y_1^2 + y_2^2
$$
in a weighted projective space $\P^{n+2}$
where the  $x_i$ have weight one and the
$y_i$ have weight $d$. {\sl Then there is an
isomorphism
$$
  H^n_o(X)\otimes T \map H^{n+2}_o(X(2)),
$$
where $T$ is a trivial Hodge structure 
of dimension one and type $(1,1)$ and
where the subscript denotes primitive cohomology}.

For the proof we note that the map
$$
    { A\Omega(x) \over P^{q+1} }
      \mapsto
    { A\Omega(x, y_1, y_2) \over ( y_1^2 + y_2^2 +P )^{q+2} }  
$$
is well-defined and via the residue provides an isomorphism
compatible with the Hodge filtrations which is defined over the
complex numbers.  However, it can be defined geometrically
and so is defined over the integers.  To see why, consider
first the trivial case $g(x) = f(x) + y_1^2 + y_2^2 = 0$ in affine
coordinates, where $x$ is a scalar variable and $f$ has degree
$2d$.  Thus $f(x) = 0$ defines a finite point set, and
$g(x) = 0$ is its double suspension.  Let $p$ be one point
of the given finite set. Then $f(p) = 0$, so the locus
$\sett{ (p,y_1,y_2) }{ y_1^2 + y_2^2 = 0 }$ lies on the double
suspension.  This locus is a pair of lines meeting in a point, and
the statement remains true in projective coordinates.  Thus we
may associate to $p$ a difference of lines $\ell_p - \ell'_p$.
This map induces an isomorphism $H_0(X,\Z) \map H_2(X(2),\Z)$
which is in fact a morphism of Hodge structures.  For the general
case we parametrize the construction just made.  The map in cohomology
which corresponds to the previous construction is the dual of the
inverse of the map in homology.

\subheading{(C) Generalizations}

The main theorem \xref{maintheorem} can be generalized in a number of
ways.  First, using the techniques of \cite{Tu},  it is certainly
possible to get  sharp results for various kinds of weighted
hypersurfaces, just as  we have obtained sharp results for standard
hypersurfaces. Second,  one can prove a quite general (but not
sharp) result that reflects the fairly weak hypothesis of criterion
\xref{largekernelcriterion}:

\proclaim{Theorem} Let $L$ be a positive line bundle on a projective
algebraic manifold $M$ of dimension at least three.  Let $P$ be the
projectivization of the  space of sections of $L^d$, and let $\Delta$
be the discrimant locus defined by sections of $L^d$ whose zero set $Z$
is singular. Then for $d$ sufficiently large the kernel of
the monodromy representation of $\Phi = \pi_1(P - \Delta)$ is large
and its image is a lattice.
\endproclaim

The monodromy representation has the primitive cohomology
$$
   H^{m-1}(Z)_0 = \kernel{\left[ H^{m-1}(Z) \mapright{Gysin} H^{m+1}(M)
\right]}
$$
as underlying vector space. The results needed for the proof are all in
the literature.  First, note that the condition that a section $s$ of
$L^d$ have a singularity of type ``$x^3 + y^3 + z^4 + \hbox{sum of
squares}$'' at a given point is  set of linear conditions and so can be
satisfied for $d$ sufficiently large.  Consequently by the
Beauville-Ebeling-Janssen argument, the image of the natural monodromy
representation is a lattice.  Second, by the results of Green
\cite{GreenPM}, the local Torelli theorem for cyclic covers holds for $d$
sufficiently large, so some component of the second monodromy
representation has nonzero differential. The standard argument used
just following Theorem
\xref{practicaldensitycriterion} proves that the discriminant locus is
irreducible, and so Theorem
\xref{udensityprop} applies to give Zariski-density for the 
second monodromy representation.  Finally, the Hodge numbers, like the
standard case of projective hypersurfaces, are polynomials in $d$ with
positive leading coefficient and of degree equal to the dimension of
$M$.  Consequently they are large for $d$ large, and therefore both the
real and complex rank of the relevant algebraic groups can be assumed
sufficiently large by taking $d$ large enough.  Thus the hypotheses of
criterion \xref{largekernelcriterion} are satisfied.

For a quick proof of the statement on the behavior of the Hodge numbers,
consider first the Poincar\'e residue sequence
$$
  0 \map \Omega^m_M \map \Omega^m_M(L^d) \map \Omega^{m-1}_Z \map 0,
$$
where $Z$ is a smooth divisor of $L^d$ and $m$ is the dimension of $M$. 
From the Kodaira vanishing theorem we have
$$
  H^0(\Omega^{m-1}_Z)_0 
    \cong 
   \cokernel{ \left[ 
      H^0(\Omega^m_M) \map H^0(\Omega^m_M(L^d)) 
        \right] }.
$$
By the Riemann-Roch theorem the dimension of the right-most term
is a polynomial with  leading coefficient  $Cd^m$, while
the dimension of the middle term is constant as a
function of $d$. Therefore the Hodge number in question
is a polynomial in $d$ of the required form.
For the other Hodge numbers we use the identification
$$
  H^q(\Omega^p_Z)_0 
    \cong
      \cokernel{ \left[ H^0(\Omega^{m-1}_Z\otimes\Theta_M\otimes
N_Z^{q-1} )
         \map
        H^0(\Omega^{m-1}_Z\otimes N_Z^q) \right] },
$$
where $N$ is the normal bundle of $Z$ in $M$, where $\Theta_M$
is the holomorphic tangent bundle of $M$, and where $p+q=m-1$.  See Proposition 6.2,
\cite{JC}, a consequence of Green's  Koszul cohomology formula (Theorem
4.f.1 in \cite{GreenKC}) for $d$ sufficiently large.  Now tensor the 
Poincar\'e residue sequence with $L^{qd}$ to get
$$
  0 
   \map 
    \Omega^m_M(L^{qd}) 
     \map 
      \Omega^m_M(L^{(q+1)d})
       \map 
        \Omega^{m-1}_Z \otimes N_Z^q
         \map
          0 .
$$
From the Kodaira vanishing theorem and the Riemann-Roch formula one
finds that the dimension of the right-hand part of the cokernel formula
is a polynomial with leading term $C((q+1)d)^m$, where
$C = c_1(L)^m/m!$.  A similar argument shows that the  dimension of
the left-hand part is a polynomial with leading coefficient
$C(qd)^m$.  Thus the leading term of $\dim H^q(\Omega^p)_0$ is bounded
below by a positive constant depending on $L$, $q$, and $M$, times $d^m$.

\subheading{(D) Questions.}

We close with some open questions.  The main
problem is, of course, to understand the nature of
the groups $\Phi_{d,n}$.  Are they linear?  
Are they residually finite?  It
seems reasonable to conjecture that in general
they are not linear groups, and, in particular,
are not lattices in Lie groups.  We settle this last question
for $\Phi_{3,2}$ in the note \cite{ACT}.

The structure of $\Phi_{d,n}$ is 
closely related to the structure of the
kernel $K$ of the natural monodromy
representation.  For $n=0$, $K$ is the pure $d$-strand
braid group the sphere and so is finitely
generated.   For $(d,n) = (3,2)$, the case of cubic
surfaces, $K$ is not finitely generated (see \cite{ACT}).  
It is therefore natural to ask when $K$ is finitely generated 
and when it is not.

\bibliography

\parskip=2pt
\baselineskip=10.5pt

\bi{ACT} D. Allcock, J. Carlson, and D. Toledo,
    Complex Hyperbolic Structures for Moduli of Cubic
    Surfaces, C. R. Acad. Sci. Paris {\bf 326}, ser I, pp 49-54 (1998)
(alg-geom/970916)

\bi{Arnold}  V.I. Arnold, S.M. Gusein-Zade, and A.N. Varchenko, Singularities of
Differentiable Maps, vol II, Birkhauser, Boston, 1988.

\bi{Beardon} A.F. Beardon, The Geometry of Discrete Groups, Springer-Verlag
(1983), pp 333.

\bi{BeauvilleLattice} A. Beauville, Le groupe de monodromie d'hypersurfaces
et d'intersections compl\`etes. Springer Lecture Notes in Mathematics,
{\bf 1194} (1986) 8--18.

\bi{BorelDT} A. Borel, Density and maximality of arithmetic groups, J. Reine 
Andgew. Math. {\bf 224} (1966) 78--89.

\bi{BHC} A. Borel and Harish-Chandra, Arithmetic subgroups of algebraic
groups, Ann. Math. (2) {\bf 75} (1962) 485--535.

\bi{CKM} H. Clemens, J. Koll\'ar, and S. Mori, editors, Higher Dimensional Complex
Geometry, Ast\'erisque {\bf 166} (1988).

\bi{JC} J. Carlson, Hypersurface Variations are Maximal, II, Trans.
Am. Math. Soc. {\bf 323} (1991), 177--196.

\bi{CG} J. Carlson and P.A. Griffiths, Infinitesimal variations of 
Hodge structure and the global Torelli problem, Journ\'ees de
G\'eometrie Alg\'ebrique D'Angers 1979, Sijthoff \& Nordhoff, Alphen
an den Rijn, the Netherlands  (1980),  51--76.

\bi{DeMo} P. Deligne, Un th\'eoreme de finitude pour la monodromie, in Discrete
Groups in Geometry and Analysis, R. Howe (ed) Birkh\"auser, Boston.

\bi{DeWeOne} P. Deligne, La Conjecture de Weil, I, Pub. Math. I.H.E.S. {\bf 43},
273 -- 307 (1974).

\bi{DeWeTwo} P. Deligne, La Conjecture de Weil, II, Pub. Math. I.H.E.S. {\bf
52}, 137 -- 252 (1980).

\bi{DeHodgeTwo} P. Deligne, Theorie de Hodge II, Publ. Math. IHES {\bf 40}
(1971), 5--58.

\bi{DOZ} G. Detloff, S. Orekov, and M. Zaidenberg, Plane curves with a
big fundamental group of the complement, Preprint 1996.

\bi{Dolgachev} I. Dolgachev, Weighted projective varieties, in Group
Actions and Vector Fields, Proceedings 1981, Lecture Notes in Math.
{\bf 956}, Springer-Verlag, New York (1982).

\bi{DolgLib} I. Dolgachev and A. Libgober, On the fundamental group to the
complement of a discriminant variety, Springer Lecture Notes in Mathematics
{\bf 862} (1981), pp 1--25.

\bi{Ebeling} W. Ebeling, An arithmetic characterisation of the symmetric
monodromy groups of singularities, Invent. Math. {\bf 77} (1984) 85-99.

\bi{Givental} A. B. Givental, Twisted Picard-Lefschetz Formulas, Funct.
Analysis Appl. {\bf 22} (1988) 10-18.

\bi{GreenKC} M. Green, Koszul cohomology and the geometry of projective
varieties. J. Differential Geom. {\bf 19} (1984), no. 1, 125--171.

\bi{GreenPM} M. Green, The period map for hypersurface sections of high
degree of an arbitrary variety. Compositio Math. {\bf 55} (1985), no. 2,
135--156.

\bi{GriffHyperbolic} P.A.  Griffiths, oral communication 1969.  

\bi{GriffPerRat} P.A. Griffiths, On the periods of certain rational 
integrals: I and II, 
Ann. of Math. {\bf 90} (1969) 460-541.

\bi{GriffPerDom} P.A. Griffiths, Periods of integrals of algebraic manifolds,
III, Pub. Math. I.H.E.S. {\bf 38} (1970), 125--180.

\bi{GS} P.A. Griffiths, and W. Schmid, Locally homogenous complex manifolds,
Acta Math. {\bf 123} (1969), 253--302.

\bi{Janssen} W. A. M. Janssen, Skew-symmetric vanishing lattices and their
monodromy groups, Math. Annalen {\bf 266} (1983) 115-133, {\bf 272} (1985) 
17-22.

\bi{Lib} A. Libgober, On the fundamental group of the space of cubic surfaces, 
Math. Zeit. {\bf 162} (1978) 63--67.

\bi{Mag} W. Magnus and A. Peluso, On a Theorem of V.I. Arnol'd, 
Comm. of Pure and App. Math. {\bf 22}, 683--692 (1969).

\bi{MargulisRigidity} G.A. Margulis, Discrete groups of motions
of manifolds of non-positive curvature, Amer. Math. Soc. Translations
{\bf 109} (1977), 33-45.

\bi{Milnor} J. Milnor, Singular Points of Complex Hypersurfaces, Annals of
Mathematics Studies {\bf 61}, Princeton University Press (1968), pp 122.

\bi{Mostow} G. D. Mostow, On a remarkable class of polyhedra in
complex hyperbolic space, Pacific Journal of Mathematics {\bf 80} (1980),
171 -- 276.

\bi{Pham} F. Pham, Formules de Picard-Lefschetz g\'en\'eralis\'ees 
et ramification des int\'egrales, Bull. Soc. Math. France 
{\bf 93}, 1965, 333-367.  Reprinted in Homology and 
Feynman Integrals, Rudolph C. Hwa and Vidgor L. Teplitz,
W. A. Benjamin 1966, pp. 331.

\bi{Seidel} P. Seidel, Floer homology and the symplectic isotopy problem,
thesis, Oxford University, 1997.

\bi{ST} M. Sebastiani and R. Thom, Un r\'esultat sur la monodromie,
Invent. Math. {\bf 13} (1971), 90--96.

\bi{SimpsonHiggs} C. T. Simpson,  Higgs bundles and local systems, Publ.
Math. IHES {\bf 75} (1992), 5--95.

\bi{Tits} J. Tits, Free subgroups in linear groups, J. Algebra {\bf 20},
250--270 (1972).

\bi{Tu} L. Tu, Macaulay's theorem and local Torelli for weighted hypersurfaces,
Compositio Math. {\bf 60} (1986), 33--44.

\bi{Zar} O. Zariski, On the problem of existence of algebraic functions
of two variables possessing a given branch curve, Amer. J. Math. {\bf 51}
(1929)  305--328.

\bi{Zimmer} R. J. Zimmer, Ergodic Theory and Semisimple Groups,
Birkh\"auser 1984, pp 209.

\endbibliography

\bigskip

\begingroup
\obeylines
\parskip=0pt
\parindent1cm
\baselineskip=10pt
\def\vs{\vskip7pt}

Department of Mathematics
University of Utah
Salt Lake City, Utah 84112
\vs

carlson@math.utah.edu $\qquad$ toledo@math.utah.edu
\vs

http://www.math.utah.edu/$\sim$carlson

http://xxx.lanl.gov --- alg-geom/9708002

To appear in Duke J. Math.

\endgroup

\enddoc
\end